\documentclass[9pt]{article}
\usepackage{latexsym,amsfonts,euscript,amsthm}

\usepackage{amssymb}
\usepackage{stmaryrd}
\usepackage{mathrsfs,amsmath}
\usepackage{leftidx}

\usepackage[figuresright]{rotating}
\usepackage[table]{xcolor}
\usepackage{multirow}
\usepackage{booktabs}
\usepackage{tabularx}

\usepackage{tocbibind} 
\usepackage[colorlinks,pagebackref]{hyperref}

\hypersetup{%
  colorlinks=true,
linktoc=all,
 citecolor=green,
    linkcolor=red,
  urlcolor=blue
}

\newtheorem{lem}{\bf Lemma}[section]

\newtheorem{prop}[lem]{\bf Proposition}

\newtheorem{thm}[lem]{\bf Theorem}
\newtheorem{rmk}[lem]{\bf Remark}
\newtheorem*{rmkA}{\bf Remark}

\newtheorem{cor}[lem]{\bf Corollary}
\newtheorem{hy}[lem]{\bf Hypothesis}

\newtheorem{mainthm}{Theorem}

\newtheorem{maincor}[mainthm]{Corollary}

\newcommand{\ackname}{Acknowledgements}
\makeatletter
\if@titlepage
  \newenvironment{acknowledgement}{%
    \titlepage
    \null\vfil
    \@beginparpenalty\@lowpenalty
    \begin{center}%
      \bfseries \ackname
      \@endparpenalty\@M
    \end{center}}%
  {\par\vfil\null\endtitlepage}
\else
  \newenvironment{acknowledgement}{%
    \if@twocolumn
      \section*{\ackname}%
    \else
      \small
      \begin{center}%
        {\bfseries \ackname\vspace{-0.5em}\vspace{\z@}}%
      \end{center}%
      \quotation
    \fi}
    {\if@twocolumn\else\endquotation\fi}
\fi
\makeatother

\makeatletter 
\@addtoreset{equation}{section}
\makeatother  

\makeatletter
\def\thanks#1{\protected@xdef\@thanks{\@thanks\protect\footnotetext{#1}}}

				   \AtEndDocument{\bigskip{\footnotesize
				   \textsc{Guohua Qian, Department of Mathematics, Suzhou University of Technology, No.99 South Third Ring
				   Road,
				   Changshu, Jiangsu, 215500, China.} \par  
				   \textit{Email address}: \href{mailto:ghqian2000@163.com}{\textsf{ghqian2000@163.com}} \par
				   \addvspace{\medskipamount}
				   \textsc{Yu Zeng, Department of Mathematics, Suzhou University of Technology, No.99 South Third Ring
				   Road,
				   Changshu, Jiangsu, 215500, China.} \par  
				   \textit{Email address}: \href{mailto:yuzeng2004@163.com}{\textsf{yuzeng2004@163.com}} 
				 }}

\title{Finite groups in which every irreducible character has either $p'$-degree or $p'$-codegree
\thanks{\textbf{Keywords}\,\, character theory, finite groups, character degrees, character codegrees.\\
\textbf{2020 MR Subject Classification}\,\, Primary 20C15, Secondary 20D20\\
Project was supported by the  NSF of China (No.12171058, 12301018), the Natural Science Foundation of Jiangsu Province  (No. BK20231356) and the Natural Science Foundation of the Jiangsu Province Higher Education Institutions of China (No. 23KJB110002).}}

\author{Guohua Qian, Yu Zeng}

\date{}

\begin{document}
\maketitle

\begin{abstract}
	For an irreducible complex character $\chi$ of a finite group $G$, the \emph{codegree} of $\chi$ is defined as $|G:\ker(\chi)|/\chi(1)$, where $\ker(\chi)$ denotes the kernel of $\chi$.
	Given a prime $p$,
	we provide a classification of finite groups
	in which every irreducible complex character has either $p'$-degree or $p'$-codegree.
\end{abstract}

\section{Introduction}

For an irreducible complex character $\chi$ of a finite group $G$,
the \emph{codegree} of $\chi$ is defined as
$$\mathrm{cod}(\chi) =\frac{|G: \ker(\chi)|}{\chi(1)}.$$
This notion was first introduced and studied in a slightly different form by D. Chillag
and M. Herzog \cite{chillag89}, and by D. Chillag, A. Mann and O. Manz \cite{chillag91}.
The current form used today was established by the first author \cite{qian02}
and was first systematically studied by the first author, Y. Wang and H. Wei \cite{qian07}.

Since the papers by I.M. Isaacs and D. Passman in the 1960s, the influence of the set of character degrees on the structure of finite groups has been extensively studied.
Many interesting results and problems have emerged from this area.
Surprisingly, some of these results and problems also have corresponding codegree versions,
leading to a wealth of interesting new theorems.
One of the main problems is the codegree analogue of Huppert's conjecture (\cite[Problem 20.79]{kourovka notebook}),
which suggests that every nonabelian finite simple group is
determined by the set of its character codegrees.
Recent papers \cite{hung25,moreto24,tongviet25} have made significant progress on this conjecture.
Furthermore,  
the set of character codegrees has been shown to have remarkable connections with element orders 
of finite groups \cite{akhlaghi24,chen22,giannelli24,isaacs11,madanha23,qian11,qian21}.

Recently, there has been a growing interest in exploring the structure of finite groups
by comparing character degrees with character codegrees.
One area of particular interest is studying the structure of finite groups $G$ 
by the set 
of the greatest common divisors of $\chi(1)$ and $\mathrm{cod}(\chi)$ for all nonlinear $\chi \in \mathrm{Irr}(G)$, denoted as:
\[
 \mathrm{GCD}(G)=\{ \mathrm{gcd}(\chi(1),\mathrm{cod}(\chi))\mid \chi \in \mathrm{Irr}(G)~\text{s.t.}~\chi(1)>1 \}.  
\]
On one hand,
finite groups $G$ where every element $\mathrm{gcd}(\chi(1), \mathrm{cod}(\chi))$ in $\mathrm{GCD}(G)$ equals $\chi(1)$ (equivalently, $\chi(1)\mid \mathrm{cod}(\chi)$ for each $\chi \in\mathrm{Irr}(G)$)
were classified by S.M. Gagola and M.L. Lewis \cite{gagola99};
for a given prime $p$, 
finite groups $G$ where the maximal $p$-power divisor of every element $\mathrm{gcd}(\chi(1), \mathrm{cod}(\chi))$ in $\mathrm{GCD}(G)$ equals
the maximal $p$-power divisor of $\chi(1)$ were characterized by the first author \cite{qian12}.
On the other hand,
finite groups $G$ where every element $\mathrm{gcd}(\chi(1), \mathrm{cod}(\chi))$ in $\mathrm{GCD}(G)$ equals $1$ (such groups are called \emph{$\mathcal{H}$-groups} in the language of \cite{liang16}) were classified by D. Liang and the first author \cite{liang16}.

Continuing this line of exploration, we study 
the $p$-analogue version of $\mathcal{H}$-groups, i.e. $\mathcal{H}_p$-groups.
Given a prime $p$,
we call a finite group $G$ an
\emph{$\mathcal{H}_p$-group} 
if 
the maximal $p$-power divisor of every element $\mathrm{gcd}(\chi(1), \mathrm{cod}(\chi))$ in $\mathrm{GCD}(G)$ equals 1, that is, every irreducible character of $G$ has either $p'$-degree or $p'$-codegree.
The study of $\mathcal{H}_p$-groups extends not only \cite[Theorem A]{qian07}, which classifies finite groups in which every nonlinear irreducible character has $p'$-codegree,
but also the celebrated It\^o-Michler theorem \cite[Theorem 5.4]{michler86}, which fully describes finite groups in which every nonlinear irreducible character has $p'$-degree.

In the next theorem, we give a complete classification of $\mathcal{H}_p$-groups.
 Before stating it, we recall two key definitions.
First, recall that a subgroup $H$ is a \emph{T.I. subgroup} (trivial intersection subgroup) of a finite group $G$ if for every $g\in G$, either $H^g=H$ or $H^g\cap H=1$.
	Second, for the special linear group $\mathrm{SL}_2(p^f)$ (where $p$ is a prime), an $\mathrm{SL}_2(p^f)$-module $V$ over the field $\mathbb{F}_p$ with $p$ elements is called the \emph{natural module} for $\mathrm{SL}_2(p^f)$
	if $V$ is isomorphic to the \emph{standard module} for $\mathrm{SL}_2(p^{f})$, i.e. the 2-dimensional vector space over the field $\mathbb{F}_{p^f}$ with $p^f$ elements (or any of its Galois conjugates) acted upon by matrix multiplication, viewed as an $\mathbb{F}_p[\mathrm{SL}_2(p^f)]$-module
	(see \cite[Definition 3.11]{parkerbook}).

\begin{mainthm}\label{thmA}
	Let $G$ be  a finite
	group and let $p$ be a prime.
	Set $N=\mathbf{O}^{p'}(G)$ and $V=\mathbf{O}_{p}(N)$.
	Then every irreducible character of $G$ has either $p'$-degree or $p'$-codegree
	if and only if one of the following holds.
	\begin{description}
		\item[(1)] $G$ has an abelian normal Sylow $p$-subgroup.
		\item[(2)]
			$N=N' \rtimes P$ where $P$ is a cyclic T.I. Sylow $p$-subgroup of $N$.
		\item[(3)] $p=3$, and $N$ is isomorphic to the affine special linear group $\mathrm{ASL}_2(3)$.
		\item[(4)] 
			$N/V$ is a Frobenius group with complement of order $p$ and cyclic kernel $K/V$
			of order $\frac{p^{pm}-1}{p^{m}-1}$,
			and $K$ is a Frobenius group with elementary abelian kernel $V$ of order $p^{pm}$.
			\item[(5)] $N$ is a nonabelian simple group, and one of the following holds.
		\begin{description}
			\item[(5a)] $p>2$, and $N$ has a cyclic Sylow $p$-subgroup.
			\item[(5b)] $N\cong \mathrm{PSL}_2(q)$ with $q=p^f$ and $f \geq 2$.
			\item[(5c)] $(N,p)\in \{ (\mathrm{PSL}_3(4),3),(M_{11},3),({}^2F_4(2)',5) \}$.
		\end{description}
		\item[(6)] $p>2$, $\mathbf{O}_{p'}(N)>1$, and $N/\mathbf{O}_{p'}(N)$ is a nonabelian simple group, and one of the following holds.
		\begin{description}
			\item[(6a)] $N$ has a cyclic T.I. Sylow $p$-subgroup. 
			\item[(6b)] $N\cong \mathrm{SL}_2(q)$ with $q=p^f$ and $f\geq 2$.
			\item[(6c)] $p=3$, and $N$ is a perfect central extension of $\mathbf{O}_{p'}(N)$ by $N/\mathbf{O}_{p'}(N)\cong \mathrm{PSL}_{3}(4)$.
		\end{description}
		\item[(7)] $V=\mathbf{C}_{N}(V)$, and one of the following holds. 
		\begin{description}
			\item[(7a)] $N/V\cong \mathrm{SL}_2(q)$ where $q=p^f\geq 4$, and $V$ is the natural module for $N/V$.
			\item[(7b)] $p=3$, $N=V \rtimes H$ where $H\cong \mathrm{SL}_2(13)$ and $V$ is a $6$-dimensional irreducible $\mathbb{F}_3[H]$-module.
			\item[(7c)] $p=3$,
			$N=V \rtimes H$ where $H\cong \mathrm{SL}_2(5)$ and $V$ is a $4$-dimensional irreducible $\mathbb{F}_3[H]$-module.	
		\end{description}
	\end{description}
\end{mainthm}

\begin{rmkA}
	{\rm We make several remarks on Theorem \ref{thmA}.
	\begin{description}
		\item[$\bullet$] If $N$ in case (2) is nonsolvable, 
		we will see in Theorem \ref{thm: description thmA (2)} that $|P|=p$.
		\item[$\bullet$] If $N$ in case (2) is solvable, 
		then $\log_p(|P|)$ cannot be bounded. 
		For instance:
		let $n\geq 2$, let $\ell=kp^n-1$ be an odd prime (the existence of $\ell$ is guaranteed by Dirichlet prime number theorem \cite[Page 169, (b)]{isaacs76}) and let $L$ be an extraspecial $\ell$-group of order $\ell^3$ with exponent $\ell$;
		as $A:=\{ a\in \mathrm{Aut}(L)\mid z^a=z~\text{for all}~z\in \mathbf{Z}(L) \}$ is isomorphic to $\mathrm{SL}_{2}(\ell)$,
		we take a Singer cycle $C(\cong \mathsf{C}_{\ell+1})$ of $A$, and take $P\leq C$ of order $p^n$;
		so $N:=L \rtimes P$ satisfies (2). 
		\item[$\bullet$] If case (4) holds, we will see in Lemma \ref{lem: 2-Frob} that $V$ is, in fact, minimal normal in $K$.
		\item[$\bullet$] If subcase (6a) holds, 
		we will see in Theorem \ref{thm: (6a) of thmA} that either $N$ is quasisimple with a cyclic Sylow $p$-subgroup,
		or $\mathbf{O}^{p'}(P\mathbf{O}_{p'}(N))$, where $P\in \mathrm{Syl}_{p}(G)$, satisfies (2). 
		 \item[$\bullet$] If subcase (7b) holds,
		 we will see in Remark \ref{rmk: smallgroup} that 
		 there are exactly two non-isomorphic groups $N$.
		 \item[$\bullet$] If subcase (7c) holds,
		 we will see in Remark \ref{rmk: smallgroup} that such $N$ is unique up to isomorphism.
	 \end{description}
	 }
\end{rmkA}

To prove Theorem \ref{thmA}, we rely on two significant prerequisites. 
The first is the deep Brauer's Height Zero Conjecture, specifically \cite[Theorem 1.1]{kessar13} and \cite[Theorem A]{malle21}. 
The second is the classification of the finite groups of order divisible by $p$ that act faithfully and irreducibly on an $\mathbb{F}_p$-module having all orbits of $p'$-size
 (\emph{$p$-exceptional linear
groups} in the language of \cite{giudici16}).
In the next section, we will present two partial results from this classification in Theorems \ref{thm: liebeck 1} and \ref{thm: liebeck 2}, which will suffice for proving the main results of this paper.

 A special class arising in the classification of $\mathcal{H}_p$-groups is 
  the class of finite groups in which every irreducible character has either $p'$-degree or $p$-defect zero.
 We call groups in this class \emph{$\mathcal{H}^*_p$-groups}.
  Y. Liu \cite{liu22} classified nonsolvable $\mathcal{H}^*_2$-groups.
  In the next corollary, we present a full classification of $\mathcal{H}^*_p$-groups.

\begin{maincor}\label{corB}
	Let $G$ be a finite group and let $p$ be a prime.  
	Set $N=\mathbf{O}^{p'}(G)$.
	Then 
	every irreducible character of $G$ has either $p'$-degree or $p$-defect zero
	if and only if one of the following holds.
	\begin{description}
		\item[(1)] $G$ has an abelian normal Sylow $p$-subgroup.
		\item[(2)] $N=N' \rtimes P$ where $P$ is a cyclic T.I. Sylow $p$-subgroup of $N$.
		\item[(3)] $N$ is a nonabelian simple group, and one of the following holds.
		\begin{description}  
			\item[(3a)] $p>2$, and $N$ has a cyclic Sylow $p$-subgroup.
			\item[(3b)] $N\cong \mathrm{PSL}_2(q)$ with $q=p^f$ and $f \geq 2$.
			\item[(3c)] $(N,p)\in \{ (\mathrm{PSL}_3(4),3),(M_{11},3),({}^2F_4(2)',5) \}$.
		\end{description}
	\item[(4)] $p>2$, $\mathbf{O}_{p'}(N)>1$, and $N/\mathbf{O}_{p'}(N)$ is a nonabelian simple group, and one of the following holds.
		\begin{description}
			\item[(4a)] $N$ has a cyclic T.I. Sylow $p$-subgroup. 
			\item[(4b)] $N\cong \mathrm{SL}_2(q)$ with $q=p^f$ and $f\geq 2$.
			\item[(4c)] $p=3$, and $N$ is a perfect central extension of $\mathbf{O}_{p'}(N)$ by $N/\mathbf{O}_{p'}(N)\cong \mathrm{PSL}_{3}(4)$.
		\end{description}
	\end{description}
\end{maincor}

Throughout this paper, we only consider finite groups and complex characters.
The paper is organized as follows: in Section 2, we gather auxiliary results; in Section 3, we collect necessary basic results on $\mathcal{H}_p$-groups; 
in Section 4, we first prove Corollary \ref{corB}
assuming Theorems \ref{thmA} and \ref{thm: simple Hp gp},
and then prove Theorem~\ref{thmA} by dealing with the $p$-solvable case in Theorem \ref{thm: classification of p-sol Hp} and the non-$p$-solvable case in Theorem \ref{thm: classification of nonsol Hp-gp} separately.

\section{Auxiliary results}

We mainly follow the notation from \cite{isaacs76} for character theory and \cite{gorenstein94} for
finite simple groups.
Throughout, we consistently refer to $p$ as a prime.
For a positive integer $n$ and a prime $p$, we write $n_p$ to
denote the maximal $p$-power divisor of $n$.
Let $G$ be a finite group. 
We use $G^\sharp$ to denote the set of nontrivial elements of $G$, $\pi(G)$ to denote the set of prime divisors of $|G|$, and $\mathrm{M}(G)$ to denote the Schur multiplier of $G$.
Let $N\unlhd G$ and $\theta \in \mathrm{Irr}(N)$.
We identify $\chi\in \mathrm{Irr}(G/N)$ with its inflation and view $\mathrm{Irr}(G/N)$ as a subset
of $\mathrm{Irr}(G)$.
We also use $\mathrm{Irr}(G|\theta)$ to denote the set of irreducible characters of $G$ lying over $\theta$,
and $\mathrm{Irr}(G|N)$ to denote the complement of the set $\mathrm{Irr}(G/N)$ in the set $\mathrm{Irr}(G)$.
Instead of $\mathrm{Irr}(G|G)$, we use $\mathrm{Irr}(G)^\sharp$ to denote the set of nontrivial irreducible characters of $G$.
Furthermore, we use $\mathsf{C}_n$ to denote a cyclic group of order $n$, 
$\mathsf{ES}(2^{1+4}_{-})$ (sometimes $2^{1+4}_{-}$) to denote the extraspecial $2$-group which is a central product of the dihedral group $\mathsf{D}_8$ and the quaternion group $\mathsf{Q}_8$,
and $\mathrm{ASL}_n(q)$ for the affine special linear group of degree $n$ over the finite field $\mathbb{F}_q$ of $q$ elements.
Other notation will be recalled or defined when necessary.

We begin by recalling some elementary results.

\begin{lem}\label{lem: order divides degree}
	Let $G$ be a finite perfect group and $\lambda\in\mathrm{Irr}(\mathbf{Z}(G))$.
	Then $o(\lambda)$, the determinantal order of $\lambda$, divides $\chi(1)$
	 for every $\chi \in \mathrm{Irr}(G|\lambda)$.
  \end{lem}
  \begin{proof}
	Let $\chi \in \mathrm{Irr}(G|\lambda)$.
	Then $\chi_{\mathbf{Z}(G)}=\chi(1)\lambda$.
	So, 
	$\det(\chi)_{\mathbf{Z}(G)}=\lambda^{\chi(1)}$.
	As $G=G'$, $\det (\chi)=1_G$.
	Therefore, $\lambda^{\chi(1)}=1_{\mathbf{Z}(G)}$, i.e. $o(\lambda)\mid \chi(1)$.
  \end{proof}

\begin{lem}\label{lem: G=opG=op'G}
	Let $G$ be a finite group with $G=\mathbf{O}^{p'}(G)$.
    Then the following hold.
\begin{description}
	\item[(1)] If $K/L$ is a $G$-chief factor of order $p$, then $K/L$ is central in $G/L$.
	\item[(2)] Given a normal series $1 \unlhd  L\unlhd  K\unlhd  G$ for $G$,
	$[G,K]\leq L$ if and only if $[P,K]\leq L$ for some $P\in \mathrm{Syl}_{p}(G)$.
\end{description}
\end{lem}
\begin{proof}
	Assume that $K/L$ is a $G$-chief factor of order $p$. 
    Set $\overline{G}=G/L$.
	As $G=\mathbf{O}^{p'}(G)$, we have $\overline{G}=\mathbf{O}^{p'}(\overline{G})$.
	 Note that $\overline{G}/\mathbf{C}_{\overline{G}}(\overline{K})$ is isomorphic to 
	 a subgroup of $\mathrm{Aut}(\overline{K})\cong \mathsf{C}_{p-1}$,
	 and so $\overline{G}=\mathbf{C}_{\overline{G}}(\overline{K})$, i.e. $\overline{K}\leq \mathbf{Z}(\overline{G})$.

	Statement (2) follows directly from the fact that $\mathbf{O}^{p'}(G)$ is the normal closure of $P \in \mathrm{Syl}_{p}(G)$ in $G$.
\end{proof}

\begin{lem}\label{lem: maximal subgroup of SL}
	Let 
	$G=\mathrm{SL}_2(p)$ where $p^2\equiv 1~(\mathrm{mod}~5)$.
	If $G$ has a subgroup $H\cong \mathrm{SL}_{2}(5)$, then $H$ is a maximal subgroup of $G$.
\end{lem}
\begin{proof}
	Assume that $G$ has a subgroup $H\cong \mathrm{SL}_{2}(5)$.
	Since $\mathrm{SL}_2(q)$ (with $q$ odd) has a unique involution $z$ 
    and $\langle z\rangle =\mathbf{Z}(\mathrm{SL}_2(q))$,
	we have $\mathbf{Z}(G)=\mathbf{Z}(H)\cong \mathsf{C}_{2}$.
    Setting $\overline{G}=G/\mathbf{Z}(G)$,
	we have $\overline{G}=\mathrm{PSL}_2(p)$, where $p^2\equiv 1~(\mathrm{mod}~5)$, and $\overline{H}\cong \mathsf{A}_5$.
	Thus, $\overline{H}$ is a maximal subgroup of $\overline{G}$ by \cite[Kapitel II, 8.27 Satz]{huppert67}.
	Consequently, $H$ is a maximal subgroup of $G$.
\end{proof}

\begin{lem}\label{lem: SL}
	Let $G$ be a finite group and let $V$ be a minimal normal subgroup of $G$.
	Assume that $G/V\cong \mathrm{SL}_{2}(q)$, where $q=p^f\geq 4$, and $|V|=q^2$.
	Then $G$ acts transitively on $V^\sharp$ 
	if and only if $V$ is the natural module for $G/V$.
\end{lem}
\begin{proof}
	Set $\overline{G}=G/V$, and let $\overline{P}$ be a Sylow $p$-subgroup of $\overline{G}$.

  If $V$ is the natural module for $\overline{G}$,
	then $G$ acts transitively on $V^\sharp$ by \cite[Lemma 3.13]{parkerbook}.

    Assume now that $G$ acts transitively on $V^\sharp$.
	Observe that $|V|=q^2<q^3$.
	By \cite[Lemma 3.12]{parkerbook}, to see that $V$ is the natural module for $\overline{G}$, it suffices to show that $|\mathbf{C}_{V}(\overline{P})|=q$.
	Since $\overline{G}\cong \mathrm{SL}_{2}(q)$ acts transitively on $V^\sharp$ and $|V|=q^2$, $\mathbf{C}_{\overline{G}}(v)\in \mathrm{Syl}_{p}(\overline{G})$ for each $v\in V^\sharp$.
    Note that $\overline{P}$ is a T.I. Sylow $p$-subgroup of $\overline{G}$, and so 
	\[
	   |V|-1= |\bigcup_{\overline{Q}\in {\rm Syl}_p(\overline{G})}
	   \mathbf{C}_{V}(\overline{Q})|-1=|{\rm Syl}_p(\overline{G})| 
	   (|\mathbf{C}_{V}(\overline{P})|-1).	
	\]
    As $|{\rm Syl}_p(\overline{G})|=q+1$, we deduce by calculation that $|\mathbf{C}_{V}(\overline{P})|=q$.
\end{proof}

\begin{lem}\label{lem: [P,G]=1}
	Let a finite group $H$ act coprimely on a finite group $G$.
	Then $H$ fixes every element of $\mathrm{Irr}(G)$ if and only if $[H,G]=1$.  
\end{lem}
\begin{proof}
	If $[H,G]=1$, then for all $h\in H$, $g\in G$ and $\chi \in \mathrm{Irr}(G)$,
	we have 
	$\chi^h(g)=\chi(g^{h^{-1}})=\chi(g)$. 
	So, $H$ fixes every element of $\mathrm{Irr}(G)$.
	
	Conversely, assume that $H$ fixes every element of $\mathrm{Irr}(G)$.
    Let $h\in H$.
    By Brauer's permutation lemma (\cite[Theorem 6.32]{isaacs76}), $h$ also fixes the conjugacy class $g^G$ for all $g\in G$.
	Given that $\mathrm{gcd}(o(h),|g^G|)=1$, 
	it follows by \cite[Lemma 13.8]{isaacs76} that
	$h$ must fix some element in $g^G$.
	Consequently, $G=\bigcup_{g\in G} \mathbf{C}_{G}(h)^g$,
	which implies that 
    $G=\mathbf{C}_{G}(h)$ for every $h\in H$.
	Therefore, $[H,G]=1$.
\end{proof}


Let $V$ be an $n$-dimensional vector space over the prime field $\mathbb{F}_p$.
As in \cite{manzwolfbook}, we denote by $\Gamma(V )$ the \emph{semilinear group} of $V$, i.e. (identifying $V$ with $\mathbb{F}_{p^n}$)
\[
	\Gamma(V)=\{ x\mapsto ax^\sigma\mid  x\in \mathbb{F}_{p^n}, a\in \mathbb{F}_{p^n}^\times, \sigma\in \mathrm{Gal}(\mathbb{F}_{p^n}/\mathbb{F}_p)\}. 
\]
It is noteworthy that $\Gamma(V)$ is a metacyclic group.

The following two theorems are partial results of the classification theorem mentioned in the introduction.
They play a crucial role in the proof of our main results.

\begin{thm}\label{thm: liebeck 1}
	Let $p$ be an odd prime and let $G$ be a finite $p$-solvable group of order divisible by $p$.
    Suppose that $V$ is a finite-dimensional, primitive $\mathbb{F}_p[G]$-module such that every orbit of $G$ on $V$ has size coprime to $p$.
	Then either $G$ is isomorphic to a subgroup of $\Gamma (V)$, or $G$ is transitive on $V^{\sharp}$. 
\end{thm}
\begin{proof}
	This is a partial result of \cite[Theorem 1]{giudici16}.
\end{proof}

\begin{thm}\label{thm: liebeck 2}
	Let $G$ be a nontrivial finite group and let $p$ be an odd prime.
	Assume that $G=\mathbf{O}^{p'}(G)=\mathbf{O}^p(G)$ and that $G$ has abelian Sylow $p$-subgroups.
	Suppose that $V$ is a finite-dimensional, primitive $\mathbb{F}_p[G]$-module such that every orbit of $G$ on $V$ has size coprime to $p$.
	Then one of the following holds.
\begin{description}
	\item[(1)] $G$ acts transitively on $V^\sharp$ and
	\begin{description}
		\item[(1a)] either $(G,|V|)=(\mathrm{SL}_2(q), q^{2})$ for some $q=p^f>4$,
		\item[(1b)] or $(G,|V|)\in \{ (2^{1+4}_{-} \cdot \mathsf{A}_5, 3^4),(\mathrm{SL}_2(13), 3^6)\}$.
	 \end{description} 
	 \item[(2)] $(G,|V|)$ is one of the following.
	 \begin{description}
		 \item[(2a)] $(G,|V|)=(\mathrm{SL}_2(5), 3^4)$ with orbit sizes $1,40,40$.
		 \item[(2b)] $(G,|V|)=(M_{11}, 3^5)$ with  orbit sizes $1,22,220$.
		 \item[(2c)] $(G,|V|)=(\mathrm{PSL}_2(11),3^{5})$ with orbit sizes $1,11,11,55,55,110$. 
\end{description}
\end{description}
\end{thm}
\begin{proof}
	This is a partial result of \cite[Theorem 5]{giudici16}.
\end{proof}

Let $G$ be a finite group 
and let $V$ be a finite dimensional $\mathbb{F}_p[G]$-module for some prime $p$.
We will assign to $V$ a finite additive group $H^{2}(G,V)$, the so-called 
\emph{second cohomology group of $V$}.
It is well-known that if $H^{2}(G,V)=0$
then every extension of $V$ by $G$ splits (see, for instance, \cite[Kapitel I, 17.2 Satz]{huppert67}).
For further details, we refer to \cite[Kapitel I, \S 16 and \S 17]{huppert67}.

\begin{rmk}\label{rmk: smallgroup}
	{\rm We supplement Theorem~\ref{thm: liebeck 2} with several observations obtained from computations in \textsf{GAP}~\cite{gap}.
	\begin{description}
		\item[$\bullet$] The group $2^{1+4}_- \cdot \mathsf{A}_5$ is the group $\mathsf{SmallGroup}(1920, 241003)$ in the $\mathsf{GAP}$ Library of small groups \cite{gap}.
		It has a unique (up to isomorphism) faithful 4-dimensional 
		irreducible module $V$ over $\mathbb{F}_3$.
		Moreover, $2^{1+4}_- \cdot \mathsf{A}_5$ acts transitively on $V^\sharp$, and $H^2(2^{1+4}_- \cdot \mathsf{A}_5,V)=0$.
		Let $\Gamma=V  \rtimes (2^{1+4}_- \cdot \mathsf{A}_5)$.
		Then there is a faithful $\chi \in \mathrm{Irr}(\Gamma)$ such that $\chi(1)=240$ and $3\mid \mathrm{cod}(\chi)$.
		Thus, $\Gamma$ is not an $\mathcal{H}_3$-group.
		\item[$\bullet$] The group $\mathrm{SL}_2(13)$ has exactly two non-isomorphic faithful 6-dimensional irreducible modules over $\mathbb{F}_3$, denoted by $V_1$ and $V_2$.
		Additionally, $\mathrm{SL}_2(13)$ acts transitively on $V_i^\sharp$, and $H^2(\mathrm{SL}_{2}(13),V_i)=0$ for $i=1,2$.
		\item[$\bullet$] The group $\mathrm{SL}_2(5)$ has a unique (up to isomorphism) faithful 4-dimensional 
		irreducible module $V$ over $\mathbb{F}_3$,
		and the orbit sizes of $\mathrm{SL}_2(5)$ on $V$ are $1,40,40$.
		Furthermore, $H^2(\mathrm{SL}_{2}(5),V)=0$.
		\item[$\bullet$] The group $M_{11}$ has exactly two non-isomorphic 5-dimensional irreducible modules over $\mathbb{F}_3$, denoted by $V_1$ and $V_2$.
		For each module $V_i$,
		$H^2(M_{11},V_i)=0$ holds.
		Let $\Gamma_i=V_i \rtimes M_{11}$.
		Then there are some $\chi_i \in \mathrm{Irr}(\Gamma_i)$
		such that $3\mid \mathrm{gcd}(\chi_i(1),\mathrm{cod}(\chi_i))$.
		Consequently, neither $\Gamma_1$ nor $\Gamma_2$ is an $\mathcal{H}_3$-group.
		\item[$\bullet$] The group $\mathrm{PSL}_2(11)$ has exactly two non-isomorphic 5-dimensional irreducible modules over $\mathbb{F}_3$, denoted by $V_1$ and $V_2$.
		For each module $V_i$,
		$H^2(\mathrm{PSL}_2(11),V_i)=0$ holds, and the orbit sizes of $\mathrm{PSL}_2(11)$ on $V_i$ are 1, 11, 11, 55, 55, 110.
		Let $\Gamma_i=V_i \rtimes \mathrm{PSL}_2(11)$.
		Then there are some faithful $\chi_i \in \mathrm{Irr}(\Gamma_i)$ such that $\chi_i(1)=33$
        and $3\mid \mathrm{cod}(\chi_i)$.
		Consequently, neither $\Gamma_1$ nor $\Gamma_2$ is an $\mathcal{H}_3$-group. 
	\end{description}
	}
\end{rmk}

We end this section with two results related to nonabelian finite simple groups.
It is noteworthy that we do not view the Tits group ${}^{2}F_4(2)'$ as a finite simple group of Lie type.

\begin{lem}[\mbox{\cite[Lemma 2.3]{qian04}}]\label{lem: out aut}
	Let $G$ be a finite almost simple group with socle $S$. 
	If $p$ divides both $|S|$ and $|G:S|$, then $G$ has
	nonabelian Sylow $p$-subgroups.
  \end{lem}

In the notation of \cite[Definition 2.2.4]{gorenstein94}, every finite simple group of Lie type is written as ${}^{d}\Sigma(q)$, where $d\in\{1,2,3\}$, $\Sigma$ is one of the root system types $A_n\;(n\ge 1)$, $B_n\;(n\ge 2)$, $C_n\;(n\ge 2)$, $D_n\;(n\ge 4)$, $E_6$, $E_7$, $E_8$, $F_4$, or $G_2$, and the parameter $q$ satisfies $q^{d}=\ell^{f}$ for some prime $\ell$ ($\ell$ is called the \emph{defining characteristic} of ${}^{d}\Sigma(q)$) and positive integer $f$.  
  On one hand, the groups $\Sigma(q):={}^{1}\Sigma(q)$ are known as the \emph{untwisted groups of Lie type}, which include: $A_n(q)$, $B_n(q)$ ($n\geq 2$), $C_n(q)$ ($n\geq 2$), $D_n(q)$ ($n\geq 4$), $E_6(q)$, $E_7(q)$, $E_8(q)$, $F_4(q)$, $G_2(q)$.
  On the other hand,
 the groups ${}^{d}\Sigma(q)$ with $d>1$ are referred to as the \emph{twisted groups of Lie type}, which include: 
 ${}^2A_n(q)$ ($n\geq 2$), ${}^2D_n(q)$ ($n\geq 4$), ${}^3D_4(q)$, ${}^2E_6(q)$, ${}^2B_2(q)$, ${}^2F_4(q)$, ${}^2G_2(q)$.
 Thus, every finite simple group of Lie type belongs to exactly one of these families.
  Also, with the exception of the Suzuki group ${}^{2}B_{2}(q)$ and the Ree groups ${}^{2}F_{4}(q)$, ${}^{2}G_{2}(q)$,
  the symbol $q$ is simply a power of the defining prime $\ell$.
  Specifically, for ${}^{2}B_{2}(q)$ and ${}^{2}F_{4}(q)$ one has $q=2^{\frac{f}{2}+1}$ (with $f$ odd), while for ${}^{2}G_{2}(q)$ one has $q=3^{\frac{f}{2}+1}$ (with $f$ odd).
  We refer to \cite{carter72,gorenstein94} for results on finite simple groups of Lie type.

  Let $S\cong{}^{d}\Sigma(q)$ be a finite simple group of Lie type in characteristic $\ell$.
  Then $S$ is defined over $\mathbb{F}_{q^d}$, where $q^d=\ell^f$.
  It is known that the automorphism group of $S$ has the structure $\mathrm{Aut}(S)=\mathrm{Inndiag}(S) \rtimes \Phi\Gamma$,
  where $\mathrm{Inndiag}(S)$ is generated by $S=\mathrm{Inn}(S)$ and the outer diagonal automorphisms of $S$,
  $\Phi\cong \mathrm{Gal}(\mathbb{F}_{q^d}/\mathbb{F}_\ell)$ induces field automorphisms, and
  $\Gamma$ induces graph automorphisms.
  For further details, we refer to \cite[Chapter 2, \S 2.5]{gorenstein94}.

\begin{lem}\label{lem: root subgp}
	Let $G$ be a finite almost simple group with socle $S$.
	Assume that $G=S \rtimes P$ where $S$ is a $p'$-group and $P$ is a cyclic Sylow $p$-subgroup of $G$.
    Then $S$ is a simple group of Lie type, and there is a nontrivial $P$-invariant abelian subgroup $A$ of $S$ such that 
	$|\mathrm{Irr}(A)|=|\mathbf{C}_{\mathrm{Irr}(A)}(Q)|^{|Q|}$ for each $Q\leq P$.  
\end{lem}
\begin{proof}
	By Feit-Thompson theorem, we know that $p>2$.
	As the outer automorphism groups of alternating groups, sporadic groups, and the Tits group ${}^2F_4(2)'$ 
	are $2$-groups (see \cite{atlas}), it follows by the classification of the finite simple groups (CFSG) that $S$ is a simple group of Lie type. 
    Let $S\cong{}^{d}\Sigma(q)$ be in defining characteristic $\ell$, where $q^d=\ell^f$.
	Write
	$\mathrm{Aut}(S)=\mathrm{Inndiag}(S) \rtimes \Phi\Gamma$
	where $\mathrm{Inndiag}(S)$, $\Phi$ and $\Gamma$ are described proceeding this lemma.
	One checks by \cite{atlas} that $\pi(\mathrm{Inndiag}(S))\cup\pi(\Gamma) \subseteq \pi(S)$.
	So, there is some $\sigma \in \mathrm{Aut}(S)$ such that $P^\sigma\leq \Phi$.
	Note that 
	$S$ has a $P$-invariant abelian subgroup $A$ such that 
	$|\mathrm{Irr}(A)|=|\mathbf{C}_{\mathrm{Irr}(A)}(Q)|^{|Q|}$ for each $Q\leq P$
	if and only if $S^{\sigma}$ has a $P^{\sigma}$-invariant abelian subgroup $A^{\sigma}$ such that 
	$|\mathrm{Irr}(A^{\sigma})|=|\mathbf{C}_{\mathrm{Irr}(A^{\sigma})}(Q^{\sigma})|^{|Q^{\sigma}|}$ for each $Q^{\sigma}\leq P^{\sigma}$.
	Without loss of generality,
	we may assume that $P\leq \Phi$.

    Let $A$ be an abelian $P$-invariant subgroup of $S$ and let $Q\leq P$.
	Note that the cyclic $p$-group $Q$ acts on $A$,
	and so Brauer's permutation lemma yields that $|\mathbf{C}_{A}(Q)|=|\mathbf{C}_{\mathrm{Irr}(A)}(Q)|$.
    Therefore, to see that $|\mathrm{Irr}(A)|=|\mathbf{C}_{\mathrm{Irr}(A)}(Q)|^{|Q|}$,
	it suffices to show that $|A|=|\mathbf{C}_{A}(Q)|^{|Q|}$.

Let $\widetilde{\Sigma}$ be a root system of $S$ and let $\hat{\Sigma}$ be the set of equivalence classes of $\widetilde{\Sigma}$ defined in \cite[Definition 2.3.1]{gorenstein94}.
By \cite[Theorem 2.4.1, Remark 2.4.2, Table 2.4]{gorenstein94}, we fix a root subgroup $X_R$ (where $R\in \hat{\Sigma}$) of $S$ and its center $A=\mathbf{Z}(X_R)$ as specified in Table \ref{tab: 1}.

	\begin{table}[!h]
		\centering
		\small
	\caption{\small Specific root subgroups of finite simple groups of Lie type and their centers}   \label{tab: 1}
	 \begin{tabular}{cccc}
	\toprule[1pt]
			Group $S$ & Type of $R$  & Root subgroup & $A=\mathbf{Z}(X_R)$\\
			 \midrule[1pt]
		   $\Sigma(q),{}^{2}D_n(q) (n\geq 4), {}^{2}E_6(q),{}^{3}D_4(q)$ & $A_1$ & $X_R=\{ x_R(t)\mid t\in \mathbb{F}_q \}$ & $A=X_R\cong \mathbb{F}_q^+$ \\
		   \midrule
		   ${}^{2}A_n(q) (n\geq 3), {}^{2}F_4(q)$ & $A_1\times A_1$ & $X_R=\{ x_R(t)\mid t\in \mathbb{F}_{q^2} \}$ & $A=X_R\cong \mathbb{F}_{q^{2}}^+$ \\
		   \midrule
		   ${}^{2}A_2(q) $ & $A_2$ & $X_R=\{ x_R(t,u)\mid t,u\in \mathbb{F}_{q^2}~\text{s.t.}~u+u^q=-tt^{q} \}$ & $A\cong \mathbb{F}_q^{+}$ \\
		   \midrule
		   ${}^{2}B_2(q) $ & $B_2$ & $X_R=\{ x_R(t,u)\mid t,u\in \mathbb{F}_{q^2} \}$ & $A\cong \mathbb{F}_{q^{2}}^+$  \\
		   \midrule
		   ${}^{2}G_2(q) $ & $G_2$ & $X_R=\{ x_R(t,u,v)\mid t,u,v\in \mathbb{F}_{q^2} \}$ & $A\cong \mathbb{F}_{q^{2}}^+$ \\
		   \bottomrule[1pt]
		\end{tabular} 

		 \smallskip

		 (Here, $\mathbb{F}^+$ denotes the additive group of the field $\mathbb{F}\in \{ \mathbb{F}_q,\mathbb{F}_{q^{2}} \}$, and it is an elementary abelian $\ell$-group of order $|\mathbb{F}|$.)
	\end{table}
	
	Let $\phi\mapsto \overline{\phi}$ be the natural isomorphism
	from $\Phi$ to $\mathrm{Gal}(\mathbb{F}_{q^d}/\mathbb{F}_\ell)$,
	and let $Q\leq P$.
	Then $\overline{Q}\leq\overline{P}\leq \mathrm{Gal}(\mathbb{F}_{q^d}/\mathbb{F}_\ell)$.
	In particular, $|\overline{P}|=|P|$ and $|\overline{Q}|=|Q|$.
	Since $\mathrm{gcd}(|\overline{P}|,d)=1$, $\overline{P}$ (and likewise $\overline{Q}$) acts faithfully via its restriction on $\mathbb{F}_q$ if $S\notin \{ {}^{2}B_2(q),{}^{2}F_4(q),{}^{2}G_2(q) \}$.
  Let $P=\langle \varphi\rangle$ and $Q=\langle \varphi_0\rangle$.
  By \cite[Theorem 2.5.1]{gorenstein94}, we have
  \begin{center}
	$x_R(t)^{\varphi}=x_R(t^{\overline{\varphi}})$ (resp. $x_R(t,u)^{\varphi}=x_R(t^{\overline{\varphi}},u^{\overline{\varphi}})$, $x_R(t,u,v)^{\varphi}=x_R(t^{\overline{\varphi}},u^{\overline{\varphi}},v^{\overline{\varphi}})$),
  \end{center}
  where $\overline{\varphi} \in \overline{P} \leq \mathrm{Gal}(\mathbb{F}_{q^{d}}/\mathbb{F}_\ell)$.
  Hence each $X_R$ is $P$-invariant, and so each $A=\mathbf{Z}(X_R)$ is also $P $-invariant.

  We first assume that $S$ is isomorphic to one of the following groups: $ \Sigma(q)$, ${}^{2}D_n(q) (n\geq 4)$, ${}^{2}E_6(q)$, ${}^{3}D_4(q)$, ${}^{2}A_n(q) (n\geq 3)$ or ${}^{2}F_4(q)$.
  Let $\mathbb{F}=\mathbb{F}_{q^2}$ if $S$ is isomorphic to either ${}^{2}A_n(q) (n\geq 3)$ or ${}^{2}F_4(q)$, otherwise let $\mathbb{F}=\mathbb{F}_q$.
  Note that $A=X_R=\{ x_R(t)\mid t\in \mathbb{F} \}\cong \mathbb{F}^+$ (the additive group of the field $\mathbb{F}$) by Table \ref{tab: 1}, and so 
  \[
    \mathbf{C}_{A}(Q)=\{ x_R(t)\mid t\in \mathbb{F}~\text{s.t.}~t^{\overline{\varphi_0}}=t \}=\{ x_R(t) \mid t\in \mathbf{C}_{\mathbb{F}}(\overline{Q}) \},
  \]
  where $\mathbf{C}_{\mathbb{F}}(\overline{Q})$ denotes the fixed field of $\overline{Q}$.
  Moreover, $|\mathbf{C}_{A}(Q)|=|\mathbf{C}_{\mathbb{F}}(\overline{Q})|$.
  By Galois theory, $[\mathbb{F}:\mathbf{C}_{\mathbb{F}}(\overline{Q})]$, the degree of 
   the field extension $\mathbb{F}/\mathbf{C}_{\mathbb{F}}(\overline{Q})$, equals $|\overline{Q}|=|Q|$.
   Therefore, $|\mathbf{C}_{A}(Q)|^{|Q|}=|\mathbf{C}_{\mathbb{F}}(\overline{Q})|^{|Q|}=|\mathbb{F}|=|A|$.

  We next assume that $S$ is isomorphic to one of the following groups: ${}^2 A_2(q)$, ${}^2 B_2(q)$ or ${}^2 G_2(q)$.

	 Assume that $S\cong{}^2 A_2(q)$.
   Then $A=\mathbf{Z}(X_R)=\{ x_R(0,u)\mid u\in \mathbb{F}_{q^2}~\text{s.t.}~u^q=-u\}$ by
   \cite[Proposition 13.6.4]{carter72}.
   Note that $A\cong \mathbb{F}_q^{+}$ by Table \ref{tab: 1},
   and so
   \[
	\mathbf{C}_{A}(Q)=\{ x_R(0,u)\mid u\in \mathbb{F}_{q^2}~\text{s.t.}~u^q=-u~\text{and}~u^{\overline{\varphi_0}}=u \}=\{ x_R(0,u)\mid u\in \mathbb{F}_{q_0^2}~\text{s.t.}~u^{q_0}=-u \}\cong \mathbb{F}_{q_0}^{+},
   \]
   where $\mathbb{F}_{q_0^2}:=\mathbf{C}_{\mathbb{F}_{q^{2}}}(\overline{Q})$.
   Moreover, $|A|=q$ and $|\mathbf{C}_{A}(Q)|=q_0$.
   Since 
   $$[\mathbb{F}_q:\mathbb{F}_{q_0}]=[\mathbb{F}_{q^{2}}:\mathbb{F}_{q_0^{2}}]=[\mathbb{F}_{q^{2}}:\mathbf{C}_{\mathbb{F}_{q^{2}}}(\overline{Q})]=|\overline{Q}|=|Q|,$$
   where the third equality holds by Galois theory,
   we conclude that $|\mathbf{C}_{A}(Q)|^{|Q|}=q_0^{|Q|}=q=|A|$.

    Assume that $S\cong{}^2B_2(q)$.
   Then $A=\mathbf{Z}(X_R)=\{ x_R(0,u) \mid u\in \mathbb{F}_{q^2}\}$ by
   \cite[Proposition  13.6.4]{carter72}.
   Note that $A\cong \mathbb{F}_{q^{2}}^{+}$ by Table \ref{tab: 1},
   and so
   \[
	\mathbf{C}_{A}(Q)=\{ x_R(0,u)\mid u\in \mathbb{F}_{q^2}~\text{s.t.}~u^{\overline{\varphi_0}}=u \}=\{ x_R(0,u)\mid u\in \mathbf{C}_{\mathbb{F}_{q^2}}(\overline{Q})\}.
   \]
   Moreover, $|A|=|\mathbb{F}_{q^{2}}|$ and $|\mathbf{C}_{A}(Q)|=|\mathbf{C}_{\mathbb{F}_{q^{2}}}(\overline{Q})|$.
   Since, by Galois theory, $[\mathbb{F}_{q^{2}}:\mathbf{C}_{\mathbb{F}_{q^{2}}}(\overline{Q})]=|\overline{Q}|=|Q|$,
   we conclude that $|\mathbf{C}_{A}(Q)|^{|Q|}=|\mathbf{C}_{\mathbb{F}_{q^{2}}}(\overline{Q})|^{|Q|}=|\mathbb{F}_{q^{2}}|=|A|$. 

   Finally, assume that $S\cong{}^2G_2(q)$.
   Then $A=\mathbf{Z}(X_R)=\{ x_R(0,0,v) \mid v\in \mathbb{F}_{q^2}\}$ by
   \cite[Proposition 13.6.4]{carter72}.
   Note that $A\cong \mathbb{F}_{q^{2}}^{+}$ by Table \ref{tab: 1},
   and so
   \[
	\mathbf{C}_{A}(Q)=\{ x_R(0,0,v)\mid v\in \mathbb{F}_{q^2}~\text{s.t.}~v^{\overline{\varphi_0}}=v \}=\{ x_R(0,0,v)\mid v\in \mathbf{C}_{\mathbb{F}_{q^2}}(\overline{Q})\}.
   \]
   Moreover, $|A|=|\mathbb{F}_{q^{2}}|$ and $|\mathbf{C}_{A}(Q)|=|\mathbf{C}_{\mathbb{F}_{q^{2}}}(\overline{Q})|$.
   Since, by Galois theory, $[\mathbb{F}_{q^{2}}:\mathbf{C}_{\mathbb{F}_{q^{2}}}(\overline{Q})]=|\overline{Q}|=|Q|$,
   we conclude that $|\mathbf{C}_{A}(Q)|^{|Q|}=|\mathbf{C}_{\mathbb{F}_{q^{2}}}(\overline{Q})|^{|Q|}=|\mathbb{F}_{q^{2}}|=|A|$.
\end{proof}

\section{Basic results on $\mathcal{H}_p$-groups}

In this section, we collect some useful results on finite $\mathcal{H}_p$-groups.
 We start by presenting some known results concerning character codegrees, which will be employed freely in the following.

\begin{lem} \label{lem: basic facts on codegree}
	Let $G$ be a finite group and let $\chi\in \mathrm{Irr}(G)$.
	\begin{description}
		\item[(1)] If $N$ is a $G$-invariant subgroup of $\ker(\chi)$,
			then the codegrees of $\chi$ in $G$ and in $G/N$ coincide.
		\item[(2)] If $M$ is a subnormal subgroup of $G$, then $\mathrm{cod}(\psi)\mid\mathrm{cod}(\chi)$ for every irreducible constituent $\psi$ of $\chi_M$.
		\item[(3)] If a prime $p$ divides $|G|$, then $p$ divides
			$\mathrm{cod}(\chi)$ for some $\chi\in \mathrm{Irr}(G)$.
			\item[(4)] If $G$ is a $p$-group and $\chi\neq 1_G$, then $p$ divides $\mathrm{cod}(\chi)$.
	\end{description}
\end{lem}
\begin{proof}
	We refer to \cite[Lemma 2.1]{liang16} for the proofs of statements (1), (2) and (3).
	For statement (4), as $\chi\neq 1_G$, 
	we have $\chi(1)<|G:\ker(\chi)|$,
	so $p$ divides $|G:\ker(\chi)|/\chi(1)=\mathrm{cod}(\chi)$ because $G$ is a $p$-group.
\end{proof}

\begin{lem}\label{lem: basic facts on Hp-group}
	Let $G$ be a finite $\mathcal{H}_p$-group and let $A$, $B$, $N$ be subgroups of $G$.
	Then the following  hold.
	\begin{description}
		\item[(1)] If $A$ is subnormal in $G$ and $B$ is normal in $A$, then $A/B$ is an $\mathcal{H}_p$-group.
		\item[(2)] If $A$ is a subnormal $p$-subgroup of $G$, then $A$ is abelian.
		\item[(3)] Let $N$ be normal in $G$, let $\theta \in \mathrm{Irr}(N)$ be of codegree divisible by $p$, and set $T=\mathrm{I}_{G}(\theta)$. Then
			\begin{description}
				\item[(3a)] $T$ contains a Sylow $p$-subgroup of $G$,
					and $\varphi$ has $p'$-degree for every $\varphi \in \mathrm{Irr}(T|\theta)$.
				\item[(3b)] if $\theta$ extends to $T$, then $G/N$ has an abelian Sylow $p$-subgroup, and $T/N$ contains a unique Sylow $p$-subgroup of $G/N$.
			\end{description}
		\item[(4)] If $G=A\times B$ with $p\mid |A|$, then $B$ has an abelian normal Sylow $p$-subgroup. 
	\end{description}
\end{lem}
\begin{proof}
	(1) Let $\theta\in \mathrm{Irr}(A/B)$ be of degree divisible by $p$, and let $\chi\in \mathrm{Irr}(G)$ be lying over $\theta$.
    Since $\theta(1)\mid \chi(1)$, it follows that 
	$p\mid \chi(1)$.
	Given that $G$ is an $\mathcal{H}_p$-group, we have $p\nmid \mathrm{cod}(\chi)$.
	By Lemma \ref{lem: basic facts on codegree} (2), 
	$\mathrm{cod}(\theta)$ divides $\mathrm{cod}(\chi)$,
	and thus, 
	we deduce that $p\nmid \mathrm{cod}(\theta)$.
	Consequently, $A/B$ is an $\mathcal{H}_p$-group.

	(2)
	We may assume that $A>1$.
 	Let $\lambda\in\mathrm{Irr}(A)^\sharp$.
	Since $A$ is a $p$-group,
	it follows that $p\mid \mathrm{cod}(\lambda)$ by Lemma \ref{lem: basic facts on codegree} (4).
	Note that $A$ is also an $\mathcal{H}_p$-group by statement (1),
	and so $\lambda(1)=1$.
	Consequently, $A$ is abelian.

	(3)
	Observe that $p$ divides the codegree of $\theta$, and thus,
	by Lemma \ref{lem: basic facts on codegree} (2), $p$ also divides the codegree of every irreducible constituent of $\theta^G$.
	As a result, all irreducible constituents of $\theta^G$ have
	$p'$-degree.
	In particular, this implies that $T$ has $p'$-index in $G$ and $P\leq T$ for some $P\in \mathrm{Syl}_{p}(G)$.
	Hence, Clifford's correspondence \cite[Theorem 6.11]{isaacs76} yields that all irreducible constituents of $\theta^T$ have
	$p'$-degree.
	Thus, (3a) holds.

	Assume further that $\theta$ extends to $T$.
	Then, by Gallagher's theorem (\cite[Corollary 6.17]{isaacs76}) and (3a), every $\alpha \in \mathrm{Irr}(T/N)$ has $p'$-degree.
	Applying It\^o-Michler theorem \cite[Theorem 5.4]{michler86} to $T/N$, we conclude that $PN/N$ is not only a Sylow 
	$p$-subgroup of $G/N$ but also an abelian normal subgroup of $T/N$.

	(4)  By Lemma \ref{lem: basic facts on codegree} (3), there exists $\alpha\in
		\mathrm{Irr}(A)$ such that $p\mid \mathrm{cod}(\alpha)$.
	As $\alpha$ extends to $G$,
	$B\cong G/A$ has an abelian normal Sylow $p$-subgroup by statement (3).
\end{proof}

\begin{lem}\label{lem: basic facts on Hp-group when OpG>1}
	Let $G$ be a finite $\mathcal{H}_p$-group with a nontrivial normal $p$-subgroup $V$.
	Then the following hold.
	\begin{description}
		\item[(1)] $G/V$ has an abelian Sylow $p$-subgroup, and
		$\mathrm{I}_{G}(\lambda)/V$ contains a unique Sylow $p$-subgroup of $G/V$ for every nontrivial $\lambda \in \mathrm{Irr}(V)$.
		\item[(2)] $[\mathbf{O}_{p'}(G),\mathbf{O}^{p'}(G)]=1$.
		\item[(3)] Assume that $P\in \mathrm{Syl}_{p}(G)$ is nonabelian. 
		Then both $V$ and $\mathrm{Irr}(V)$ are irreducible $\mathbb{F}_p[G]$-modules and 
		$$|\mathrm{Irr}(V)|-1=|\mathrm{Syl}_{p}(G)|(|\mathbf{C}_{\mathrm{Irr}(V)}(P)|-1).$$
		Furthermore, if $p>2$, then both $V$ and $\mathrm{Irr}(V)$ are primitive $\mathbf{O}^{p'}(G/\mathbf{C}_{G}(V))$-modules over $\mathbb{F}_p$.
	\end{description}
\end{lem}
\begin{proof}
	(1) Let $\lambda\in \mathrm{Irr}(V)^\sharp$ and $T=\mathrm{I}_{G}(\lambda)$.
	As $V$ is a nontrivial normal $p$-subgroup of an $\mathcal{H}_p$-group $G$, it follows that $p\mid \mathrm{cod}(\lambda)$ by Lemma \ref{lem: basic facts on codegree} (4) and that $\lambda(1)=1$ by Lemma \ref{lem: basic facts on Hp-group} (2).
	So, by Lemma \ref{lem: basic facts on Hp-group} (3), in order to establish statement (1), it suffices to show that $\lambda$ extends to $T$.
	In fact, Lemma \ref{lem: basic facts on Hp-group} (3) asserts that $T$ contains a Sylow $p$-subgroup $P$ of $G$, and every $\varphi\in \mathrm{Irr}(T|\lambda)$ has $p'$-degree; 
	so $\varphi_P$ has a linear irreducible constituent $\mu$ lying over $\lambda$;
	this implies that $\lambda$ extends to $P$, and consequently, $\lambda$ extends to $Q$ for every Sylow subgroup $Q/V$ of $T/V$;
	therefore, $\lambda$ extends to $T$ by \cite[Corollary 11.31]{isaacs76}.

	(2) Set $W=\mathbf{O}_{p'}(G)$ and $K=\mathbf{O}_{p}(G)$.
	Let $\alpha\in \mathrm{Irr}(W)$ and fix a $\beta\in \mathrm{Irr}(K)^\sharp$.
	Then $\alpha\times \beta\in \mathrm{Irr}(W\times K)$ has codegree divisible by $p$ by Lemma \ref{lem: basic facts on codegree}.
	Applying Lemma \ref{lem: basic facts on Hp-group} (3),
	we have that $P\leq \mathrm{I}_{G}(\alpha\times \beta)$ for some $P\in \mathrm{Syl}_{p}(G)$.
	By statement (1), $\mathrm{I}_{G}(\beta)$ contains a unique Sylow $p$-subgroup of $G$, which we denote by $P_\beta$.
	Since $P\leq \mathrm{I}_{G}(\alpha\times \beta)=\mathrm{I}_{G}(\alpha) \cap  \mathrm{I}_{G}(\beta)$,
	it follows that $P_\beta=P\leq \mathrm{I}_{G}(\alpha)$ for each $\alpha \in \mathrm{Irr}(W)$.
    As $P_\beta$ acts coprimely on $W$,
	$[W,P_\beta]=1$ by Lemma \ref{lem: [P,G]=1}.
	So, by Lemma \ref{lem: G=opG=op'G}, we conclude that $[W,\mathbf{O}^{p'}(G)]=1$.

	(3) Assume that $P\in \mathrm{Syl}_{p}(G)$ is nonabelian.
	Set $U=\mathrm{Irr}(V)$ and $C=\mathbf{C}_{G}(U)$.
	By Lemma \ref{lem: basic facts on Hp-group} (2), $V$ and hence $U$ are abelian $p$-groups.
	Consider now the action of $G$ on $U$.
	We assert that $P$ acts nontrivially on $U$.
	In fact, otherwise $\mathrm{I}_{G}(\lambda)$ contains every Sylow $p$-subgroup of $G$;
	so, statement (1) forces that $P\unlhd G$;
	however, $P$ is nonabelian which contradicts Lemma \ref{lem: basic facts on Hp-group} (2).
	Note that $PC/C$ acts faithfully, nontrivially on $U$
	and that $\mathrm{I}_{G}(\lambda)/C$ contains a unique Sylow $p$-subgroup of $G/C$ for each $\lambda \in U^\sharp$ by statement (1).	 
	Therefore, \cite[Lemma 4]{zhang00} implies that $U$ is an irreducible $\mathbb{F}_p[G]$-module.
	 By \cite[Lemma 1]{zhang00}, $V$ is also an irreducible $\mathbb{F}_p[G]$-module and $\mathbf{C}_{G}(V)=C$.
	Recalling that every nontrivial element of $U$ 
	 is fixed by a
	 unique Sylow $p$-subgroup of $G$ by statement (1), we deduce that $U=\bigcup_{Q\in {\rm Syl}_p(G)} \mathbf{C}_{U}(Q)$, and
	 $\mathbf{C}_{U}(Q_1)\cap \mathbf{C}_{U}(Q_2)=\{ 1_V \}$ whenever $Q_1, Q_2$ are distinct Sylow $p$-subgroups of $G$. 
	 By calculation,
	 $$|U|-1=|\bigcup_{Q\in {\rm Syl}_p(G)}
	 \mathbf{C}_{U}(Q)|-1=|{\rm Syl}_p(G)| 
	 (|\mathbf{C}_{U}(P)|-1).$$ 
	 Next, we assume that $p>2$.
	 Set $X/C=\mathbf{O}^{p'}(G/C)$.
	 As $PC/C\in \mathrm{Syl}_{p}(X/C)$, $\mathrm{I}_{X}(\lambda)/C$ also contains a unique Sylow $p$-subgroup of $G/C$ for each $\lambda \in U^\sharp$.
	 Again, by \cite[Lemma 4]{zhang00}, $U$ is a primitive $\mathbb{F}_p[X/C]$-module.
	 Consequently, $V$ is also a primitive $\mathbb{F}_p[X/C]$-module by \cite[Lemma 1]{zhang00}.
\end{proof}

Let $G$ be a finite $p$-solvable group and let $l_p(G)$ denote its \emph{$p$-length}.
It is well-known that $l_p(G/\mathbf{O}_{p',p}(G))$ equals $l_p(G)-1$ when $p\mid |G|$,
and that $l_p(G/N)=l_p(G)$ for any normal subgroup $N$ of $G$ contained in $\Phi(G)$ or $\mathbf{O}_{p'}(G)$.
For further details, we refer to \cite[Kapitel VI, \S 6]{huppert67}.

\begin{lem}\label{lem: basic facts on Hp p-sol gp}
	Let $G$ be a finite $p$-solvable $\mathcal{H}_p$-group.
	Then the following hold.
	\begin{description}
		\item[(1)] $l_p(G)\leq 2$.
		\item[(2)] If $l_p(G)\leq 1$, then $G$ has an abelian Sylow $p$-subgroup.
		\item[(3)] If $\mathbf{O}_{p}(G)>1$, then either $P\in \mathrm{Syl}_{p}(G)$ is abelian, or
			$\mathbf{O}_p(G)$ is minimal normal in $G$ and $\mathbf{O}_p(G)\cap
				\Phi(G)=1$.
	\end{description}
\end{lem}
\begin{proof}
	(1) 
	By induction, we may assume $\mathbf{O}_{p'}(G)=1$ and $\mathbf{O}_{p}(G)>1$.
	Hence, it follows by Lemma \ref{lem: basic facts on Hp-group} (2)
	and Lemma \ref{lem: basic facts on Hp-group when OpG>1} (1) that $\mathbf{O}_{p}(G)$ is abelian, and $G/\mathbf{O}_{p}(G)$ has an abelian Sylow $p$-subgroup.
	Consequently, $l_p(G)\leq 2$.

	(2) Since $l_p(G)\leq 1$, there exist $G$-invariant subgroups $N\leq M$
	such that $M/N$
	is isomorphic to a Sylow $p$-subgroup of $G$.
Therefore, the desired result follows directly from Lemma \ref{lem: basic facts on Hp-group}.

	(3) 
	Assume that $\mathbf{O}_{p}(G)>1$ and that $P \in \mathrm{Syl}_{p}(G)$ is nonabelian. 
	By Lemma \ref{lem: basic facts on Hp-group when OpG>1}, 
	$[\mathbf{O}_{p'}(G),\mathbf{O}^{p'}(G)]=1$ and,
	$\mathbf{O}_{p}(G)$ is minimal normal in $G$.
	Thus, $\mathbf{O}_{p',p}(G)=\mathbf{O}_{p'}(G)\times \mathbf{O}_{p}(G)$ and,
	either $\mathbf{O}_{p}(G)\cap \Phi(G)=1$ or $\mathbf{O}_{p}(G)\leq \Phi(G)$.
    If $\mathbf{O}_{p}(G)\leq \Phi(G)$, then $l_p(G)=l_p(G/\mathbf{O}_{p}(G))=l_p(G/\mathbf{O}_{p',p}(G))\leq 1$
	where the two equalities hold by \cite[Kapitel VI, 6.4 Hilfssatz]{huppert67} and the inequality holds by statement (1),
	whereas statement (2) implies that $P$ is abelian, a contradiction.
	Therefore, $\mathbf{O}_{p}(G)\cap \Phi(G)=1$.
\end{proof}

We end this section with some facts on $\mathcal{H}^*_p$-groups.
Before that, we briefly introduce some facts on blocks.

Let $p$ be a prime and let $G$ be a finite group.
Then $\mathrm{Irr}(G)$ is a disjoint union of $\mathrm{Irr}(B)$ (the \emph{set of irreducible characters in $B$}) with $B$ running over all $p$-blocks of $G$.

Suppose that $B$ is a $p$-block of $G$ with defect group $D$.
R. Brauer proved that $|G:D|_p$ is the maximal power of $p$ dividing the degrees of all 
characters in $\mathrm{Irr}(B)$.
So, there is some $\chi \in \mathrm{Irr}(B)$ with $\chi(1)_p=|G:D|_p$.
If $D=1$, then we call the block $B$ has \emph{defect zero} (in this case, $B$ contains exactly one irreducible character $\chi$ which has $p$-defect zero in $G$);
if $D$ is a Sylow $p$-subgroup of $G$, then we call the block $B$ has \emph{maximal defect}.
Even though $G$ may not have a defect zero $p$-block,
it always has a maximal defect $p$-block.
For instance, the \emph{principal $p$-block} of $G$ (denoted by $B_0$), which is
the unique $p$-block containing the principal character of $G$, always has maximal defect.

A celebrated result of J.A. Green states that every defect group of a $p$-block of $G$ is an intersection of two Sylow $p$-subgroups of $G$ (see, for instance, \cite[Corollary 4.21]{navarrobook}).
So, 
\begin{center}
	 	$G$ has a T.I. Sylow $p$-subgroup $\Rightarrow$ every $p$-block of $G$ has either maximal defect or defect zero.
\end{center}
 In general, the converse of the above statement does not hold.
 However, if we assume that a Sylow $p$-subgroup of $G$ is abelian, the converse is true
 (see \cite[Theorem 3.2]{pazderski91}).
 This leads us to the following lemma.

\begin{lem}\label{lem: max defect or defect zero}
	Suppose that a finite group $G$ has an abelian Sylow $p$-subgroup $P$.
	Then  $P$ is a T.I. subgroup of $G$ if and only if every $p$-block of $G$ has either maximal defect or defect zero.
\end{lem}

Applying Lemma \ref{lem: max defect or defect zero}, \cite[Theorem 1.1]{kessar13} 
and \cite[Theorem A]{malle21}, we obtain a rough characterization of $\mathcal{H}^*_p$-groups via 
their Sylow $p$-subgroups.
Recall that an $\mathcal{H}^*_p$-group is a finite group in which every irreducible character has either $p'$-degree or $p$-defect zero.
It is also important to note that $\mathcal{H}^*_p$-groups are indeed $\mathcal{H}_p$-groups.

\begin{prop}\label{prop: tisylow}
	Let $G$ be a finite group.
	Then $G$ has an abelian T.I. Sylow $p$-subgroup if and only if 
	$G$ is an $\mathcal{H}^*_p$-group.
  \end{prop}
\begin{proof}
	Assume that $G$ has an abelian T.I. Sylow $p$-subgroup.
	Then, by Lemma \ref{lem: max defect or defect zero}, every $p$-block of $G$ has either maximal defect or defect zero. 
	Let $B$ be a $p$-block of $G$ with maximal defect.
	By \cite[Theorem 1.1]{kessar13}, every $\chi\in \mathrm{Irr}(B)$ has $p'$-degree,
	which implies that $G$ is an $\mathcal{H}^*_p$-group.
	
	Now, assume that $G$ is an $\mathcal{H}^*_p$-group.
	Then every $p$-block of $G$ has either maximal defect or defect zero. 
	Let $B_0$ be the principal $p$-block of $G$.
	Since every $\chi \in \mathrm{Irr}(B_0)$ has $p'$-degree,
	\cite[Theorem A]{malle21} implies that $G$ has an abelian Sylow $p$-subgroup. 
	Therefore, by Lemma \ref{lem: max defect or defect zero}, $G$ has an abelian T.I. Sylow $p$-subgroup.
\end{proof}

\begin{cor}\label{cor: subgroup is Hp*}
	Let $G$ be a finite group and let $P\in \mathrm{Syl}_{p}(G)$. Then the following hold.
	\begin{description}
		\item[(1)] Assume that $P\leq H\leq G$.
		If $G$ is an $\mathcal{H}^*_p$-group,
		then $H/N$ is an $\mathcal{H}^*_p$-group 
		whenever $N\unlhd H$. 
		\item[(2)] Let $Z$ be a central $p'$-subgroup of $G$. If $G/Z$ is an $\mathcal{H}^*_p$-group, then $G$ is also an $\mathcal{H}^*_p$-group.
		\item[(3)] $G$ is an $\mathcal{H}^*_p$-group if and only if $\mathbf{O}^{p'}(G)$ is an $\mathcal{H}^*_p$-group.
	\end{description}
\end{cor}
\begin{proof}
	 Note that $P$ is abelian implies that $PN/N$ is abelian for $N\unlhd G$,
	  and that $PN/N$ is abelian for a normal $p'$-subgroup $N$ of $G$ implies that $P$ is abelian. 
	  So, in the proofs of statements (1) and (2), we only need to verify the T.I. property of Sylow $p$-subgroups.

	(1) Let $N\unlhd H$.
	  Then $PN/N$ is a Sylow $p$-subgroup of $H/N$.
	  For each $x\in H$, as $G$ is an $\mathcal{H}^*_p$-group,
	  Proposition \ref{prop: tisylow} implies that 
	  either $P=P^x$ or $P\cap P^{x}=1$;
	  consequently, it is straightforward to verify that either 
	  $PN/N=P^xN/N$ or $PN/N\cap P^{x}N/N=1$.
	  By Proposition \ref{prop: tisylow} again, $H/N$ is also an $\mathcal{H}^*_p$-group.

	  (2) Let $x\in G$. 
	  Since $G/Z$ is an $\mathcal{H}^*_p$-group, we have $P^xZ\cap PZ=Z$ or $P^xZ=PZ$.
      If $P^xZ\cap PZ=Z$, then, as $Z$ is a $p'$-subgroup of $G$, it follows that $P^x\cap P\leq P\cap Z=1$.
	  On the other hand, if $P^xZ=PZ$, then, given that $Z\leq \mathbf{Z}(G)$, we conclude that $P^x=P$.
	  Hence, $G$ is also an $\mathcal{H}^*_p$-group.

	  (3) Assume that $N:=\mathbf{O}^{p'}(G)$ is an $\mathcal{H}^*_p$-group.
	  Let $\chi\in \mathrm{Irr}(G)$ be of degree divisible by $p$ and let $\theta$ be an irreducible constituent of $\chi_{N}$.
	  Then $\chi(1)_p=\theta(1)_p$.
	  Since $p\mid \theta(1)$, we have $\theta(1)_p=|N|_p=|G|_p$.
	  Consequently, $\chi(1)_p=|G|_p$, implying that $G$ is an $\mathcal{H}^*_p$-group. 
      The converse statement follows directly from statement (1).
\end{proof}

\section{Main results}

Note that every finite group is either $p$-solvable or non-$p$-solvable.
So, 
we split the proof of Theorem \ref{thmA} into two parts:
in Theorem \ref{thm: classification of p-sol Hp}, we consider the $p$-solvable case, while the non-$p$-solvable case will be dealt with in Theorem \ref{thm: classification of nonsol Hp-gp}.
Before that, we prove Corollary \ref{corB} assuming Theorems \ref{thmA} and \ref{thm: simple Hp gp}.

\begin{proof}[Proof of Corollary \ref{corB}]
 We first assume that $G$ is an $\mathcal{H}_p^*$-group.
	Then $G$ is an $\mathcal{H}_p$-group
	such that 
	either every character in $\mathrm{Irr}(G)$ has $p'$-degree or,
	there exists a character in $\mathrm{Irr}(G)$ having $p$-defect 0.
	If the former holds, then $G$ has an abelian normal Sylow $p$-subgroup by It\^o-Michler theorem.
	Assume that the latter holds.
	Then $\mathbf{O}_{p}(N)\leq \mathbf{O}_{p}(G)=1$.
	Applying Theorem \ref{thmA} to $G$ and omitting
	the cases with $\mathbf{O}_{p}(N)>1$ in Theorem \ref{thmA},
	we are done.
	
    Conversely, we assume that one of the cases (1), (2), (3) or (4) holds.
    If one of (1), (2) or (4a) holds,
    then $N$ has an abelian T.I. Sylow $p$-subgroup,
	so $G$ is an $\mathcal{H}_p^*$-group by Proposition \ref{prop: tisylow} and Corollary \ref{cor: subgroup is Hp*}.
	If one of (3), (4b) or (4c) holds,
	then $\mathbf{O}_{p'}(N)\leq \mathbf{Z}(N)$ and, 
	$N/\mathbf{O}_{p'}(N)$ is an $\mathcal{H}_p^*$-group by Theorem \ref{thm: simple Hp gp}.
	So, $G$ is an $\mathcal{H}_p^*$-group by Corollary \ref{cor: subgroup is Hp*}.
\end{proof}

\subsection{$p$-solvable $\mathcal{H}_p$-groups}

Given a prime $p$,
the aim of this subsection is to classify finite $p$-solvable $\mathcal{H}_p$-groups.
Note that finite $p$-solvable $\mathcal{H}_p$-groups have $p$-length
at most $2$ by Lemma \ref{lem: basic facts on Hp p-sol gp} (1).
So, our strategy is to classify them according to their $p$-length.

\begin{lem} \label{lem: Hp p-sol gp with p-length 1}
	Let $G$ be a finite $p$-solvable $\mathcal{H}_p$-group with $p$-length $1$ and let $N=\mathbf{O}^{p'}(G)$.
	Assume that $P\in \mathrm{Syl}_{p}(G)$ is not normal in $G$.
	Then the following hold.
	\begin{description}
		\item[(1)] $\mathbf{O}^{p}(N)<D$ 
		for every normal subgroup $D$ of $G$ of order divisible by $p$.
			In particular, $N'=\mathbf{O}^{p}(N)$.
		\item[(2)] If $\chi\in \mathrm{Irr}(G)$ has degree divisible by $p$, then
			$\chi$ has $p$-defect zero in $G$.  
		\item[(3)] $P$ is a cyclic T.I. subgroup of $N$. 	
	\end{description}
\end{lem}
\begin{proof}
	As $G$ is a $p$-solvable $\mathcal{H}_p$-group with $p$-length $1$,
	$N=K \rtimes P$ where $K=\mathbf{O}_{p'}(N)=\mathbf{O}^{p}(N)=\mathbf{O}^{p',p}(G)$, and $P$ is abelian by Lemma \ref{lem: basic facts on Hp p-sol gp} (2).

	(1) Let $D$ be a normal subgroup of $G$ of order divisible by $p$, and let $A/B$ be a $G$-chief factor of order divisible by $p$ within $D$. 
	Set $\overline{G}=G/B$.
	Since $\overline{G}$ is a $p$-solvable $\mathcal{H}_p$-group
	with $\mathbf{O}_p(\overline{G})>1$, Lemma \ref{lem: basic facts on Hp-group when OpG>1} (2)
	implies that $\overline{P}$ centralizes $\overline{K}$.
	So, $\overline{P}\unlhd \overline{N}$ and hence $\overline{P}\unlhd \overline{G}$.
	This means that $\mathbf{O}^{p',p}(\overline{G})=1$, so 
	$B\geq \mathbf{O}^{p',p}(G)=K$.
	Consequently, $D> K$.
	
	Observe that $N$ is also a $p$-solvable $\mathcal{H}_p$-group with $p$-length $1$, and that $P$ is not normal 
	in $N$.
	As $PN' \unlhd N$ has order divisible by $p$, it follows that $K\leq PN'$.
	Given that $N'\leq K$, we conclude that $N'=K$.

	(2) Let $\chi\in \mathrm{Irr}(G)$ be of degree divisible by $p$, and let $\theta$ be an irreducible constituent
	of $\chi_N$.
	Since $p$ does not divide $|G:N|$,
	it follows that $\theta(1)_p=\chi(1)_p$.
	Given that $N/K\cong P$ is abelian,
  every character in $\mathrm{Irr}(N/K)$ has $p'$-degree.
	Therefore, $\ker(\theta)$ does not contain $K$.
	Note that $N$ satisfies the hypotheses of this lemma,
	and so statement (1) yields that $p\nmid |\ker(\theta)|$.
	As $p\nmid \mathrm{gcd}(\theta(1),\mathrm{cod}(\theta))$, it follows that $\theta(1)_p=|P|$.
    So, $\chi(1)_p=|P|$.

	(3) By statement (2) and Proposition \ref{prop: tisylow}, $P$ is an abelian T.I. Sylow $p$-subgroup of the $\mathcal{H}_p^*$-group $G$.
	As $N= N'  \rtimes P$
	where $P$ acts nontrivially on $N'$ by statement (1),
    there is a $P$-invariant Sylow $q$-subgroup $Q$ of $N'$ such that $[Q,P]>1$.
	Let $H=\mathbf{O}^{p'}(QP)$.
	Then $H=Q_0 \rtimes P$ where $Q_0\in \mathrm{Syl}_{q}(H)$,
	and $H$ is a solvable $\mathcal{H}^*_p$-group with $p$-length 1 by Corollary \ref{cor: subgroup is Hp*} (1).
	Note that $P$ is not normal in $H$, and so $Q_0=H'$ by statement (1).
	Let $H'/E$ be an $H$-chief factor.
	According to statement (1),
	$H'/E$ is the unique minimal normal subgroup of the solvable group $H/E$.
	Applying \cite[Lemma 12.3]{isaacs76} to $H/E$, we conclude
	that $PE/E$ is cyclic.
	Therefore, $P$ must also be cyclic, as $P\cong PE/E$.
\end{proof}

\begin{lem}\label{lem: implies solvable}
	Let $G$ be a finite $p$-solvable $\mathcal{H}_p$-group with $p$-length $2$.
    Assume that $G=V \rtimes D$ where $D=\mathbf{O}^{p'}(D)$ and $V=\mathbf{O}_{p}(G)$ is the unique minimal normal subgroup of $G$.
	Then $G$ is solvable.
\end{lem}
\begin{proof}
    As a $2$-solvable group is also solvable by Feit-Thompson theorem, we may assume that $p>2$.
	Observe that $V=\mathbf{O}_{p}(G)$ is a normal $p$-subgroup of $G$.
	To see that $G$ is solvable, it remains to show that $D$ is solvable.
    
	Now, note that $l_p(G)=2$, and hence Sylow $p$-subgroups of $G$ are nonabelian. 
    Since $D=\mathbf{O}^{p'}(D)(\cong G/\mathbf{O}_{p}(G))$ is a $p$-solvable $\mathcal{H}_p$-group with $p$-length $1$ and $\mathbf{O}_{p}(D)=1$,
	it follows by Lemma \ref{lem: Hp p-sol gp with p-length 1} that
	$D=D' \rtimes P$, where $P$ is a cyclic T.I. Sylow $p$-subgroup of $D$. 
	
	Recall that Sylow $p$-subgroups of $G$ are nonabelian and that $p>2$.
    By Lemma \ref{lem: basic facts on Hp-group when OpG>1}, $D$ acts primitively on $U:=\mathrm{Irr}(V)$, and $\mathrm{I}_{D}(\lambda)$ contains a unique Sylow $p$-subgroup of $D$ for every $\lambda \in U^\sharp$. 
    Note that $D$ is $p$-solvable,
	and hence an application of Theorem \ref{thm: liebeck 1} yields that 
	either $D$ is isomorphic to a subgroup of $\Gamma(U)$ or $D$ acts transitively on $U^\sharp$.
	If the former holds, then we are done.
    
	Thus, we may assume that $D$ acts transitively on $U^\sharp$, and that $D$ is
	isomorphic to neither a subgroup of $\Gamma(U)$
	nor $\mathrm{SL}_{2}(3)$.
	So, using the classification of the 2-transitive affine permutation groups (see \cite[Appendix 1, Hering's theorem]{liebeck87}),
	we deduce that $D$ belongs to either the Extraspecial classes or the Exceptional classes (in the language of \cite[Appendix 1, Hering's theorem]{liebeck87}).
	If $D$ belongs to the Exceptional classes, then, given that $D$ is $p$-solvable, $K\unlhd D<\mathrm{GL}_2(p)$ and
	$\mathrm{SL}_2(5)\cong K\leq D\cap \mathrm{SL}_2(p)<\mathrm{SL}_2(p)$ where $p\in \{11,19,29,59\}$.
    Since $p^2\equiv 1~(\mathrm{mod}~5)$, $K$ is a maximal subgroup of $\mathrm{SL}_{2}(p)$ by Lemma \ref{lem: maximal subgroup of SL}.
    Consequently, $D\cap \mathrm{SL}_2(p)=K$, and so
	$p$ does not divide $|\mathrm{SL}_2(p)D:\mathrm{SL}_2(p)|\cdot|K|=|D:K|\cdot|K|=|D|$, a contradiction.

	On the other hand, we assume that $D$ belongs to the Extraspecial classes.
	Then $D\leq \mathbf{N}_{\mathrm{GL}(U)}(R)$ and
	either $(R,|U|)= (\mathsf{Q}_8,p^2) $ where $p\in \{ 5,7,11,23 \}$,
	or $(R,|U|)= (\mathsf{ES}(2^{1+4}_{-}),3^4)$ and $D/R\leq \mathsf{S}_5$.
	If the former holds, as $\mathrm{Aut}(\mathsf{Q}_8)\cong \mathsf{S}_4$, then $RP=R  \times P$ where $R\cong \mathsf{Q}_8$.
	However, $\mathbf{N}_{\mathrm{GL}(U)}(P)/P\cong \mathsf{C}_{p-1}\times \mathsf{C}_{p-1}$, a contradiction. 
    If the latter holds, then 
	$D/R=D'/R \rtimes PR/R$ where $|PR/R|=3$.
	So, $D'/R$ is isomorphic to neither $\mathsf{A}_5$ nor $\mathsf{S}_5$.
     Therefore, $D$ is solvable, and we are done.
\end{proof}

\begin{lem} \label{lem: Hp p-sol gp with p-length 2}
	Let $G$ be a finite $p$-solvable $\mathcal{H}_p$-group with $p$-length $2$.
	If $G=\mathbf{O}^{p'}(G)$,
	then one of the following holds.
	\begin{description}
		\item[(1)] $p=3$, and $G\cong\mathrm{ASL}_2(3)$.
		\item[(2)] $G$ has a normal series $1 \lhd V \lhd K \lhd G$
		such that $G/V$ is a Frobenius group with complement of order $p$ and cyclic kernel $K/V$
				of order $\frac{p^{pm}-1}{p^m-1}$,
				and $K$ is a Frobenius group with elementary abelian kernel $V$ of order $p^{pm}$.	
	\end{description}
\end{lem}
\begin{proof}
	Let $N=\mathbf{O}_{p'}(G)$ and $V=\mathbf{O}_{p',p}(G)$,
	and set $\overline{G}=G/N$.
	Note that $l_p(\overline{G})=2$, and hence Sylow $p$-subgroups of $\overline{G}$ are nonabelian.	
	Given that $\overline{G}$ is a $p$-solvable $\mathcal{H}_p$-group with $\mathbf{O}_{p'}(\overline{G})=1$ and $\overline{V}=\mathbf{O}_{p}(\overline{G})>1$, 
	it follows by Lemma \ref{lem: basic facts on Hp p-sol gp} (3)
	that $\overline{V}$ is the unique minimal normal subgroup of $\overline{G}$ and $\overline{V}\cap \Phi(\overline{G})=1$.
	Consequently, $\overline{G}=\overline{V}  \rtimes \overline{D}$ where $\overline{D}$ is a complement of $\overline{V}$ in $\overline{G}$.
	Furthermore, $\overline{V}$ is not cyclic,
	as this would imply $l_p(\overline{G})=1$, contradicting 
	the fact that $l_p(\overline{G})=2$. 
	As $G=\mathbf{O}^{p'}(G)$ and $\overline{D}\cong G/\mathbf{O}_{p',p}(G)$ has $p$-length 1,
	$\overline{D}=\mathbf{O}^{p'}(\overline{D})$ is a $p$-solvable $\mathcal{H}_p$-group such that $l_p(\overline{D})=1$ and $\mathbf{O}_{p}(\overline{D})=1$.
	So, Lemma \ref{lem: Hp p-sol gp with p-length 1} implies 
	that $\overline{D}=\overline{H}  \rtimes \overline{P}$ where $\overline{P}$ is a cyclic T.I. Sylow $p$-subgroup of $\overline{D}$ and $\overline{H}=\overline{D}'$.
	Moreover, by Lemma \ref{lem: implies solvable}, $\overline{G}$ is solvable.
	Set $U=\mathrm{Irr}(\overline{V})$.

	\emph{Claim 1.} Either $p=3$, $\overline{G}\cong \mathrm{ASL}_2(3)$ and $\overline{H}\cong \mathsf{Q}_8$, or $\overline{D}\leq \Gamma(U)$ and $\overline{H}$ is cyclic.

	Let $\lambda\in U^\sharp$.
	An application of Lemma \ref{lem: basic facts on Hp-group when OpG>1} (1) yields that $\mathrm{I}_{\overline{D}}(\lambda)$ contains a unique Sylow $p$-subgroup of $\overline{D}$ which is abelian.
	As $\overline{V}$ is the unique minimal normal subgroup of $\overline{G}$,
    both $\overline{V}$ and $U$ are faithful irreducible $\overline{D}$-modules by \cite[Lemma 1]{zhang00}.
	Noting that $\overline{G}$ is solvable and 
	applying \cite[Main Lemma]{palfy01} to $\overline{D}$ and $U$, we conclude that either
	$p=3$, $|U|=3^2$ and $\mathrm{SL}_2(3)\leq \overline{D}\leq \mathrm{GL}_2(3)$,
	or
	$\overline{D}\leq
	\Gamma(U)$.
	If the former holds, as $\overline{D}=\mathbf{O}^{3'}(\overline{D})$, then $\overline{D}=\mathrm{SL}_2(3)$ and $\overline{H}\cong \mathsf{Q}_8$.
	Since, up to isomorphism, $\overline{D}=\mathrm{SL}_2(3)$ has a unique $2$-dimensional irreducible module over $\mathbb{F}_3$, i.e. the natural module for $\overline{D}$ over $\mathbb{F}_3$ (check via $\mathsf{GAP}$ \cite{gap}),
	$\overline{V}$ is also isomorphic to the natural module.
    Therefore, $\overline{G}\cong\mathrm{ASL}_2(3)$.
	If the latter holds, as $\Gamma(U)'$ is cyclic, then $\overline{H}=\overline{D}'$ is cyclic.

	\emph{Claim 2.} $N=1$.

	 As $V=\mathbf{O}_{p',p}(G)$, $V= N \rtimes Q$, where $Q\in \mathrm{Syl}_{p}(V)$, is a $p$-solvable $\mathcal{H}_p$-group with $l_p(V)=1$.
	 Recall that $Q\cong \overline{V}$ is not cyclic.
	 Consequently, by Lemma \ref{lem: Hp p-sol gp with p-length 1}, $Q\unlhd V$, 
	 indicating that $Q$ is a nontrivial normal $p$-subgroup of $G$.
	 Now, applying Lemma \ref{lem: basic facts on Hp-group when OpG>1} (2) to $G$, we deduce that $N$ is central in $G=\mathbf{O}^{p'}(G)$.
	 Also, since $G=\mathbf{O}^{p'}(G)$, the $G$-invariant $p'$-subgroup $N\leq G'$.
	 As a consequence, $N$ is isomorphic to a quotient group of $\mathrm{M}(\overline{G})$ (the Schur multiplier of $\overline{G}$).
	 Since $\overline{H}\in \mathrm{Hall}_{p'}(\overline{G})$ is either a quaternion $2$-group or cyclic by Claim 1,
	 every Sylow subgroup of $\overline{H}$ has a trivial Schur multiplier. 
	 Thus $\mathrm{M}(\overline{G})$ is a $p$-group by \cite[Kapitel V, 25.1 Satz]{huppert67}, forcing $N$ to be trivial.

	\emph{Claim 3.}
	If  $D\leq \Gamma(U)$, then $G/V$ is a Frobenius group with cyclic complement $PV/V$ and cyclic kernel $HV/V$,
	and $HV$ is a Frobenius group with cyclic complement $H$ and elementary abelian kernel $V$.

	Recall that $D=H \rtimes P$ is a solvable $\mathcal{H}_p$-group with $p$-length 1 where $P\in \mathrm{Syl}_{p}(D)$ and $H$ are both cyclic and that $\mathbf{O}_{p}(D)=1$.
	Also, $V$ is an elementary abelian $p$-group.
	
	Now, consider the coprime action of $P$ on $H$.
	Let $P_0$ be a nontrivial subgroup of $P$.
	As $H=\mathbf{C}_{H}(P_0)\times [H,P_0]$, the group $P_0[H,P_0]$ is a $D$-invariant subgroup of order divisible by $p$.
	According to Lemma \ref{lem: Hp p-sol gp with p-length 1}, we have $H=\mathbf{O}^{p',p}(D)< P_0[H,P_0]$.
	It follows that $\mathbf{C}_{H}(P_0)=1$ for every nontrivial subgroup $P_0$ of the cyclic $p$-group $P$.
	Therefore, $D$ is a Frobenius group with cyclic complement $P$ and cyclic kernel $H$.
	Since $G/V\cong D$, $G/V$ is a Frobenius group with cyclic complement $PV/V$ and cyclic kernel $HV/V$.

	Next, we claim that $HV$ is a Frobenius group with cyclic complement $H$ and elementary abelian kernel $V$.
	Indeed, for
	each $\lambda\in U^\sharp$, since $\mathrm{I}_{D}(\lambda)$ contains a unique Sylow
	$p$-subgroup of $D$ 
	by Lemma \ref{lem: basic facts on Hp-group when OpG>1} (1),
	and $D$ is a Frobenius group with complement $P\in \mathrm{Syl}_{p}(D)$,
	it forces $\mathrm{I}_{D}(\lambda)\in \mathrm{Syl}_{p}(D)$;
	so $\mathrm{I}_{H}(\lambda)=1$, indicating $HV$ is a Frobenius group with cyclic complement $H$ and elementary abelian kernel $V$.

	\emph{Claim 4.}
	If  $D\leq \Gamma(U)$, then $|G/HV|=p$, $|HV/V|=\frac{p^{pm}-1}{p^m-1}$ and $|V|=p^{pm}$ where $p^m=|\mathbf{C}_{U}(P)|$.

	Let $P_0$ be a maximal subgroup of $P$.
	Note that $D$ is a Frobenius group with cyclic complement $P$ and cyclic kernel $H$ and that $\mathbf{C}_{U}(H)=\{ 1_V \}$ by Claim 3.
	According to \cite[Theorem 15.16]{isaacs76}, we have
	$$\dim_{\mathbb{F}_p}
		\mathbf{C}_{U}(P_0)=|P:P_0|\dim_{\mathbb{F}_p}
		\mathbf{C}_{U}(P)=p\dim_{\mathbb{F}_p}
		\mathbf{C}_{U}(P).$$
	Now, take $\lambda\in
		\mathbf{C}_{U}(P_0)-
		\mathbf{C}_{U}(P)$.
	 By Lemma \ref{lem: basic facts on Hp-group when OpG>1} (1),
    $D$ has a unique Sylow $p$-subgroup $Q$
	(distinct from $P$) that fixes $\lambda$.
	Since 
    $D$ is a Frobenius group with complement $P$,
	we have $P_0\leq P\cap Q=1$, implying $P_0=1$.
	Therefore, $|G/HV|=|P|=p$ and $|V|=|U|=p^{pm}$ for some positive integer $m$, where $p^{m}=|\mathbf{C}_{U}(P)|$.

	Recall that Sylow $p$-subgroups of $G$ are nonabelian, and that $G=V \rtimes D$ where $V=\mathbf{O}_{p}(G)$.
	By Lemma \ref{lem: basic facts on Hp-group when OpG>1} (3), we deduce that
	$$p^{pm}-1=|U|-1=|{\rm Syl}_p(G)|
	(|\mathbf{C}_{U}(P)|-1)=|{\rm Syl}_p(D)|
		(|\mathbf{C}_{U}(P)|-1)=|H|(p^m-1).$$
		So, $|HV/V|=|H|=\frac{p^{pm}-1}{p^m-1}$.
\end{proof}

We will see in the next lemma that the subgroup $V$ appearing in statement (2) of Lemma~\ref{lem: Hp p-sol gp with p-length 2} is, in fact, a minimal normal subgroup of the group $K$ mentioned in that statement.

\begin{lem}\label{lem: 2-Frob}
	Let $p,r$ be primes.
	Let $G$ be a finite group having a normal series $1 \lhd V \lhd K \lhd G$.
	Assume that $G/V$ is a Frobenius group with complement of order $p$ and cyclic kernel $K/V$
	of order $l=\frac{r^{pm}-1}{r^{m}-1}$,
	and that $K$ is a Frobenius group with elementary abelian kernel $V$ of order $r^{pm}$.
	Then $V$ is minimal normal in $K$.
\end{lem}
\begin{proof} 
	Let $L$ be a Frobenius complement of $V$ in $K$.
	Note that $L$ is a Hall $r'$-subgroup of the solvable group $K$,
	and hence the Frattini's argument yields that $G=K\mathbf{N}_{G}(L)=V\mathbf{N}_{G}(L)$ where $V\cap \mathbf{N}_{G}(L)=\mathbf{C}_{V}(L)=1$.
    Set $H=\mathbf{N}_{G}(L)$.
	As $H\cong G/V$, $H$ is a Frobenius group with complement $P\cong \mathsf{C}_{p}$ and kernel $L\cong \mathsf{C}_{l}$.

	Assume that $V$ is not minimal normal in $K$.
    In other words, $V_L$ (the restriction of $H$-module $V$ to $L$) is not irreducible as an $L$-module.
    Observe that $|H:L|=p$ is a prime and that $\mathbf{C}_{H}(L)=L$,
	and so \cite[Theorem 0.1, Lemma 2.2]{manzwolfbook} yields that $V_L=V_1\oplus\cdots \oplus V_p$
    where $V_i$ are irreducible $L$-modules and $H/L$ acts transitively on $\{ V_1,\dots ,V_p \}$.
	In particular, all $\mathbf{C}_{L}(V_i)$ share the same order.
	As $L$ is a cyclic group of order $l$, it follows that $\mathbf{C}_{L}(V_1)=\cdots =\mathbf{C}_{L}(V_p)$.
	Since $\mathbf{C}_{L}(V)=\bigcap_{i=1}^{p}\mathbf{C}_{L}(V_i)=1$,
    each $V_i$ is a faithful irreducible $L$-module of the cyclic group $L$.
    So, the dimension of each $V_i$ is equal to the order of $r$ modulo $l$.
	Recall that $V_L=V_1\oplus\cdots \oplus V_p$ has dimension $pm$, 
	and hence the dimension of each $V_i$ is equal to $m$.
	Therefore, $l\mid r^{m}-1$.
    By calculation,
	\[
	l= \frac{r^{pm}-1}{r^{m}-1}=r^{(p-1)m}+\cdots +r^{m}+1 \equiv p~(\mathrm{mod}~l),
	\]
	which contradicts $\operatorname{gcd}(p,l)=1$.
\end{proof}

Now, we are ready to classify finite $p$-solvable $\mathcal{H}_p$-groups.

\begin{thm}\label{thm: classification of p-sol Hp}
	Let $G$ be a finite $p$-solvable
	group and let $N=\mathbf{O}^{p'}(G)$. 
	Set $V=\mathbf{O}_{p}(N)$.
	Then
	$G$ is an $\mathcal{H}_p$-group
	if and only if one of the following holds.
	\begin{description}
		\item[(1)] $G$ has an abelian normal Sylow $p$-subgroup.
		\item[(2)]
			$N=N' \rtimes P$ where $P$ is a cyclic T.I. Sylow $p$-subgroup of $N$.
		\item[(3)] $p=3$, and $N$ is isomorphic to the affine special linear group $\mathrm{ASL}_2(3)$.
		\item[(4)] 
			$N/V$ is a Frobenius group with complement of order $p$ and cyclic kernel $K/V$
			of order $\frac{p^{pm}-1}{p^m-1}$,
			and $K$ is a Frobenius group with elementary abelian kernel $V$ of order $p^{pm}$.
	\end{description}
\end{thm}
\begin{proof}
	Let $P$ be a Sylow $p$-subgroup of $G$.
	We first assume that $G$ is an $\mathcal{H}_p$-group.
  If $P\unlhd G$, then $P$ is abelian by Lemma \ref{lem: basic facts on Hp-group} (2).
  Assume that $P$ is not normal in $G$.	
  Note that $N=\mathbf{O}^{p'}(N)$ is also a $p$-solvable $\mathcal{H}_p$-group.
  Therefore, the desired results follow directly by Lemmas \ref{lem: Hp p-sol gp with p-length 1} and \ref{lem: Hp p-sol gp with p-length 2}.

	Conversely, suppose that one of the cases (1), (2), (3) or (4) holds.
	Let $\chi\in \mathrm{Irr}(G)$ be of degree divisible by $p$, and let $\theta$ be an irreducible constituent of $\chi_N$.
    To see that $G$ is an $\mathcal{H}_p$-group, it suffices to show that $\mathrm{cod}(\chi)_p=1$.
	If (1) holds, then $G$ is an $\mathcal{H}_p$-group by It\^{o}-Michler theorem.
	If (2) holds, then $G$ is an $\mathcal{H}_p$-group by Proposition \ref{prop: tisylow}.
	If (3) holds, then $p=3$, $\theta$ is the unique irreducible character in $\mathrm{Irr}(N)$ of degree $3$ and $V\leq \ker(\theta)$ (check via $\mathsf{GAP}$ \cite{gap}).
	Note that $V\leq\ker(\theta)\leq \ker(\chi)$,
	and so
	$$\mathrm{cod}(\chi)_3=\frac{|G:\ker(\chi)|_3}{\chi(1)_3}=\frac{3}{3}=1.$$
	Assume now that (4) holds.
	We claim that $V\leq \ker(\theta)$. 
	Let $\lambda\in \mathrm{Irr}(V)^\sharp$.
	Then 
	$\mathrm{I}_{N}(\lambda)/V\in \mathrm{Syl}_{p}(N/V)$ has order $p$.
   In fact, 
	as $N/V$ is a Frobenius group with complement of order $p$ and kernel $K/V$ acting faithfully on $\mathrm{Irr}(V)$, 
	\cite[Theorem 15.16]{isaacs76} implies that $|\mathbf{C}_{\mathrm{Irr}(V)}(Q/V)|=p^{m}$ for every $Q\in \mathrm{Syl}_{p}(N)$;
	for distinct $Q_1,Q_2\in \mathrm{Syl}_{p}(N)$, note that
	$$\mathbf{C}_{\mathrm{Irr}(V)}(Q_1/V)\cap \mathbf{C}_{\mathrm{Irr}(V)}(Q_2/V)=\mathbf{C}_{\mathrm{Irr}(V)}(\langle Q_1/V,Q_2/V\rangle)\leq \mathbf{C}_{\mathrm{Irr}(V)}(K_0/V)=\{ 1_V \}$$
    where the last equality holds because $K_0$, the preimage of $K_0/V:=K/V\cap\langle Q_1/V,Q_2/V\rangle$ in $K$,
	is a Frobenius group with kernel $V$;
	therefore, $\mathrm{Irr}(V)^\sharp=\bigcup_{Q\in \mathrm{Syl}_{p}(N)}\mathbf{C}_{\mathrm{Irr}(V)}(Q/V)^\sharp$
    by comparing the sizes of the two sets;
	consequently, $\mathrm{I}_{N}(\lambda)/V\in \mathrm{Syl}_{p}(N/V)$ has order $p$.
    So, $\lambda$ extends to $\mathrm{I}_{N}(\lambda)$ by \cite[Corollary 11.22]{isaacs76}.
	Clifford's theorem and Gallagher's theorem then force that every character in $\mathrm{Irr}(N|\lambda)$
	has $p'$-degree.
	Thus, $V\leq \ker(\theta)$.
	Recall that $N/V$ is a Frobenius group with complement of order $p$ and kernel $K/V$,
	and so $\ker(\theta)<K$ and $\theta(1)=p$.
	Finally, since $V\leq \ker(\chi)\cap N\leq \ker(\theta)<K<N=\mathbf{O}^{p'}(G)$, we have 
	$$\mathrm{cod}(\chi)_p=\frac{|G:\ker(\chi)|_p}{\chi(1)_p}=\frac{|N\ker(\chi):\ker(\chi)|_p}{\theta(1)}=\frac{|N:\ker(\chi)\cap N|_p}{\theta(1)}=\frac{|N:K|}{\theta(1)}=\frac{p}{p}=1.$$
\end{proof}

Finally, we describe the groups that arise in the case (2) of Theorem \ref{thm: classification of p-sol Hp}.

\begin{lem}\label{lem: Hp* inertia}
	Let $G$ be a finite $\mathcal{H}^*_p$-group and $N=\mathbf{O}_{p'}(G)$.
	If $G/N$ is a cyclic $p$-group,
	then either $\mathrm{I}_{G}(\theta)=N$ or $\mathrm{I}_{G}(\theta)=G$ for each $\theta\in \mathrm{Irr}(N)$.
\end{lem}
\begin{proof}
	Let $\theta\in \mathrm{Irr}(N)$ be not $G$-invariant and $T=\mathrm{I}_{G}(\theta)$.
	Then $\theta(1)_p=1$ because $N=\mathbf{O}_{p'}(G)$.
    As $T/N$ is cyclic, $\theta$ extends to $\hat{\theta}\in \mathrm{Irr}(T)$,
	and $\chi:=\hat{\theta}^G\in \mathrm{Irr}(G)$ by Clifford's correspondence.
    Noting that $\chi(1)_p=|G:T|$ is a power of $p$, 
	and that $G$ is an $\mathcal{H}^*_p$-group, we conclude that $T=N$.
\end{proof}


\begin{thm}\label{thm: description thmA (2)}
	Let $G$ be a finite group.
	Assume that $G=\mathbf{O}^{p'}(G)=G' \rtimes P$ where $P$ is a cyclic T.I. Sylow $p$-subgroup of $G$.
    If $G'$ has a nonabelian $G$-chief factor $N/K$,
	then $P$ acts nontrivially on $N/K$ and 
	$|P|=p$. 
\end{thm}
\begin{proof}
    Note that, as $G=\mathbf{O}^{p'}(G)$, $P$ acts nontrivially on $N/K$ by Lemma \ref{lem: G=opG=op'G}.
		For a minimal $P$-invariant quotient group $N/L$ of $N/K$,
	as $G$ is an $\mathcal{H}^*_p$-group by Proposition \ref{prop: tisylow},	
	$\mathbf{O}^{p'}(PN/L)$ satisfies the hypotheses of this theorem by Corollary \ref{cor: subgroup is Hp*}.
    Assume that $|P|=p^n>p$, and let $G$ be a counterexample of minimal order.
    By the minimality of $G$,
	$G=N \rtimes P$ where $N=G'$ is a nonabelian minimal normal subgroup of $G$. 
    So, $N=S_1\times \cdots \times S_t$
	where $S_i$ are isomorphic to a nonabelian simple group $S:=S_1$.
	Let $\{ y_1=1,y_2,\cdots ,y_t \}~(\subseteq P)$ be a transversal of $\mathbf{N}_{G}(S)$ in $G$ and let $P=\langle y\rangle$. 
	Then $o(y)=p^{n}$.

	Assume first that $t>1$.
	Let $\alpha$ be a nontrivial $\mathrm{Aut}(S)$-invariant character in $\mathrm{Irr}(S)$ (see, for instance, \cite[Lemma 2.11]{moreto05}) 
	and let $\theta=\alpha\times (1_S)^{y_2}\times \cdots \times (1_S)^{y_t} \in \mathrm{Irr}(N)$.
	Then $\mathrm{I}_{G}(\theta)=\mathbf{N}_{G}(S)$.
	Given that $G$ is an $\mathcal{H}^*_p$-group,
	it follows by Lemma \ref{lem: Hp* inertia} that $\mathbf{N}_{G}(S)=N$.
	Thus, $\{ 1,y,y^2,\cdots ,y^{p^n-1} \}$ becomes a transversal of $\mathbf{N}_{G}(S)$ in $G$. 
	Let 
	$$\varphi=\alpha\times (1_S)^{y}\times \cdots \times (1_S)^{y^{p-1}}\times \alpha^{y^p}\times (1_S)^{y^{p+1}}\times \cdots \times (1_S)^{y^{2p-1}}\times \cdots \times  \alpha^{y^{p^n-p}}\times (1_S)^{y^{p^n-p+1}}\times \cdots \times (1_S)^{y^{p^n-1}}\in \mathrm{Irr}(N).$$
    Then $N<\mathrm{I}_{G}(\varphi)=N\langle y^p\rangle<G$ which contradicts Lemma \ref{lem: Hp* inertia}.

	Assume next that $t=1$.
	In this case, $G=N \rtimes P$ is an almost simple group with socle $N$ where $P\in \mathrm{Syl}_{p}(G)$ is a cyclic group of order $p^n$.
	Let $1<P_0<P$.
	According to Lemma \ref{lem: root subgp}, 
	$N$ has an abelian $P$-invariant subgroup $A$ 
	such that $\mathbf{C}_{\mathrm{Irr}(A)}(P)<\mathbf{C}_{\mathrm{Irr}(A)}(P_0)$.
	Set $H=AP$.
	So, we have $A<\mathrm{I}_{H}(\lambda)<H$ for each $\lambda\in \mathbf{C}_{\mathrm{Irr}(A)}(P_0)-\mathbf{C}_{\mathrm{Irr}(A)}(P)$.
	However, as $H$ is also an $\mathcal{H}^*_p$-group by Corollary \ref{cor: subgroup is Hp*} (1), we conclude a contradiction by Lemma \ref{lem: Hp* inertia}.
\end{proof}

\subsection{Non-$p$-solvable $\mathcal{H}_p$-groups}

In this subsection, we provide a classification of finite non-$p$-solvable $\mathcal{H}_p$-groups.
We begin with a theorem concerning nonabelian finite simple $\mathcal{H}_p$-groups.
It is noteworthy that, for a nonabelian finite simple group $G$,
it is an $\mathcal{H}_p$-group if and only if it is an $\mathcal{H}^*_p$-group.

\begin{thm}\label{thm: simple Hp gp}
	Let $G$ be a nonabelian finite simple group.
	Then $G$ is an $\mathcal{H}_p$-group if and only if one of the following holds.
	\begin{description}
	  \item[(1)] $p>2$, and $G$ has a cyclic Sylow $p$-subgroup.
	  \item[(2)] $G\cong \mathrm{PSL}_2(q)$ where $q=p^f$ and $f \geq 2$.
	  \item[(3)] $(G,p)\in \{ (\mathrm{PSL}_3(4),3),(M_{11},3),({}^2F_4(2)',5) \}$.
	\end{description} 
  \end{thm}
  \begin{proof}
	Assume first that $G$ is an $\mathcal{H}_p$-group.  
    Since $G$ is nonabelian simple, it must also be an $\mathcal{H}^*_p$-group.
    Applying Proposition \ref{prop: tisylow},
	we conclude that 
    $G$ has an abelian T.I. Sylow $p$-subgroup, say $P$.
	If $P$ is noncyclic, then the classification of nonabelian simple groups 
	with a noncyclic T.I. Sylow $p$-subgroup (\cite[Proposition 1.3]{blau90})
	implies that either 
	\[
	 (G,p)\in \{ (\operatorname{PSU}_3(p^{n}),p),({}^2B_2(2^{m+\frac{1}{2}}),2),({}^2G_2(3^{m+\frac{1}{2}}),3),(McL,5),(J_4,11)\},
	\]
	or one of (2) or (3) holds.
	The former case is ruled out as Sylow $p$-subgroups of $G$ are nonabelian. 
	In fact, if $(G,p)\in \{ (\operatorname{PSU}_3(p^{n}),p),({}^2B_2(2^{m+\frac{1}{2}}),2),({}^2G_2(3^{m+\frac{1}{2}}),3)\}$,
	this can be verified  
	using \cite[Proposition 13.6.4]{carter72};
	if $(G,p)\in \{(McL,5),(J_4,11)\} $, 
	this can be confirmed by referring to \cite{atlas}.

	Conversely, let us assume that one of the cases (1), (2) or (3) holds.
	If (3) holds, then $G$ is an $\mathcal{H}_p$-group by checking \cite{atlas}.
	If (1) holds, i.e. $P$ is cyclic, then \cite[Corollary 2]{blau85} yields that $G$ is an $\mathcal{H}_p$-group.
    If (2) holds, as $G$ has an abelian T.I. Sylow $p$-subgroup in this case,
	then $G$ is an $\mathcal{H}_p$-group by Proposition \ref{prop: tisylow}.
	\end{proof}

\begin{lem}\label{lem: nonsol Hp gp with nonabel minimal normal}
	Let $G$ be a finite $\mathcal{H}_p$-group with a nonabelian minimal normal subgroup $N$.
	Assume that $p$ divides $|N|$.
	Then $N=\mathbf{O}^{p'}(G)$ is a nonabelian simple $\mathcal{H}_p$-group.
\end{lem}
\begin{proof}
	Since $N$ is a nonabelian minimal normal subgroup of $G$, $N=S\times T$, where $S$ is a nonabelian simple
	group.   
	Write $C=\mathbf{C}_{G}(N)$.
	Observe that $N C=S\times TC$ is an $\mathcal{H}_p$-group,
	and hence $p\nmid |TC|$ by Lemma \ref{lem: basic facts on Hp-group} (4).
	So, $T=1$, implying that $N=S$ and $p\nmid |C|$.
    In summary, $N$ is a nonabelian simple $\mathcal{H}_p$-group of order divisible by $p$ and $|G/CN|_p=|G/N|_p$.
	
	Assume that $p\mid |G/N|$, and let $G$ be a counterexample of minimal order.
	Since $|G/CN|_p=|G/N|_p$, 
	the minimality of $G$ implies that $G$ is an almost simple group with socle $N$.
	Note that $N$ is a nonabelian simple $\mathcal{H}_p$-group of order divisible by $p$, and that $p\mid |\mathrm{Out}(N)|$.
	Therefore, 
	 by Theorem \ref{thm: simple Hp gp} and \cite{atlas}, 
	  either $N$ has a cyclic Sylow $p$-subgroup with $p>2$ or 
	 $(N,p)\in \{ (\mathrm{PSL}_2(p^f),p),(\mathrm{PSL}_3(4),3)\}$.
	 As the outer automorphism groups of alternating groups, sporadic groups, and the Tits group ${}^{2}F_4(2)'$ 
	 are $2$-groups, it follows by the CFSG that $N$ must be a simple group of Lie type.

	We claim that $\mathbf{O}_{p}(G/N)=1$. 
	Assume not.
	 Then $G$ has a subnormal subgroup $H$ such that $|H/N|=p$.
	 Note that $H$ is also an almost simple $\mathcal{H}_p$-group with socle $N$,
	 and so $G=H$ by the minimality of $G$ whence $|G/N|= p$.
	 Let $\chi\in \mathrm{Irr}(G|N)$ be of degree divisible by $p$.
	 As $N$ is the unique minimal normal subgroup of $G$, it follows that $\ker(\chi)=1$.
	 Since $p\nmid \mathrm{cod}(\chi)$, $\chi$ must have $p$-defect zero in $G$.
	 Note also that $G/N$ is abelian, and so every character in $\mathrm{Irr}(G)$ has either $p'$-degree or $p$-defect zero. 
	 Now, applying Proposition \ref{prop: tisylow} to $G$, we conclude that $G$ has an abelian Sylow $p$-subgroup which is contrary to Lemma \ref{lem: out aut}.

	 Given that $p\mid |G/N|$ and $\mathbf{O}_{p}(G/N)=1$, It\^o-Michler theorem yields  
	 the existence of some $\alpha \in \mathrm{Irr}(G/N)$
	 of degree divisible by $p$.
	  Let $\theta$ be the Steinberg character of $N$.
	    Then $\theta$ extends to $\chi \in \mathrm{Irr}(G)$ (see \cite{schmid85}).
	  By Gallagher's theorem, $\chi \alpha \in \mathrm{Irr}(G|\theta)$.
	  Given that $N$ is the unique minimal normal subgroup of $G$,
	  it follows that $\ker(\chi)=\ker(\chi\alpha)=1$.
	  Now, we conclude the following two statements:
    if $p\mid \theta(1)$, then $p\mid (\chi(1),\mathrm{cod}(\chi))$;
	  if $p\nmid \theta(1)$, then $p\mid (\chi\alpha(1),\mathrm{cod}(\chi\alpha))$.
	  In each case, we conclude a contradiction. 
	  Thus $N=\mathbf{O}^{p'}(G)$.  
\end{proof}

From Lemma \ref{lem: under hyp perfect basics} to Lemma \ref{lem: under hyp perfect N p C>1}, we will address the most complicated cases that arise in the process of
classifying finite non-$p$-solvable $\mathcal{H}_p$-groups. In order to avoid repetitions, we introduce the following hypothesis.

\begin{hy}\label{hyp: perfect}
	Let $G$ be a finite $\mathcal{H}_p$-group and let $R$ be the $p$-solvable radical of $G$ $($i.e. the maximal
	$p$-solvable normal subgroup$)$.
	Assume that $G=\mathbf{O}^{p}(G)=\mathbf{O}^{p'}(G)$ and that $G/R$ is a nonabelian simple group with $R>1$.
\end{hy}

Assuming Hypothesis \ref{hyp: perfect},
we see that $R$ is the unique maximal normal subgroup of the perfect group $G$ and that $G/R$ is a nonabelian simple $\mathcal{H}_p$-group of order divisible by $p$.

\begin{lem}\label{lem: under hyp perfect basics}
	Assume Hypothesis \ref{hyp: perfect}.
	Then $R$ does not have a $G$-chief factor of order $p$.
\end{lem}
\begin{proof}
Let $G$ be a counterexample of minimal order.
   Then $G$ has a unique minimal normal subgroup $V$, and $|V|=p$.
  Also, by Lemma \ref{lem: G=opG=op'G}, $V\leq \mathbf{Z}(G)$.
   Now, let $\chi\in \mathrm{Irr}(G|V)$ and let $\lambda$ be an irreducible constituent of $\chi_V$.
   As $G=\mathbf{O}^{p}(G)=\mathbf{O}^{p'}(G)$,
   $G$ is a perfect group.
   So, Lemma \ref{lem: order divides degree} implies that $p=o(\lambda)$ divides $\chi(1)$.
   Since $V$ is the unique minimal normal subgroup of $G$, it follows that $\ker(\chi)=1$.
   Note that $p\nmid\mathrm{gcd}(\chi(1),\mathrm{cod}(\chi))$,
   and so $\chi$ has $p$-defect zero in $G$ whereas $\chi(1)\mid |G/V|$, a contradiction.
\end{proof}

\begin{lem}\label{lem: under hyp perfect N p'}
	Assume Hypothesis \ref{hyp: perfect} and that  $p\nmid |R|$.
	Then $G$ is an $\mathcal{H}^*_p$-group,
	and
	one of the following holds.
	\begin{description}
		\item[(1)] $G$ has a cyclic T.I. Sylow $p$-subgroup.
		\item[(2)] $p>2$, $G\cong \mathrm{SL}_2(p^f)$ with $f\geq 2$.
		\item[(3)] $p=3$, $G$ is a perfect central extension of $R$ by $G/R\cong \mathrm{PSL}_{3}(4)$.
	\end{description}
\end{lem}
\begin{proof}
	Let $\chi \in \mathrm{Irr}(G|R)$ be of degree divisible by $p$.
   Observe that $R$ is the unique maximal normal subgroup of $G$,
   and so $\ker(\chi)< R$.    
   Since $p\nmid\mathrm{gcd}(\chi(1),\mathrm{cod}(\chi))$,
   we see that $\chi(1)_p=|G/\ker(\chi)|_p=|G|_p$.
   As $G/R$ is a nonabelian simple $\mathcal{H}_p$-group,
   every character in $\mathrm{Irr}(G/R)$ has either $p'$-degree or $p$-defect zero,
   so does every character in $\mathrm{Irr}(G)$.
   In other words, $G$ is an $\mathcal{H}^*_p$-group.
   Hence, Proposition \ref{prop: tisylow} guarantees the existence of an abelian T.I. Sylow $p$-subgroup of $G$.
   Let $P\in \mathrm{Syl}_{p}(G)$ and set $H=RP$.
   Then $H$ is a $p$-solvable $\mathcal{H}_p$-group with $p$-length 1 by Corollary \ref{cor: subgroup is Hp*} (1).

   Assume now that $P$ is not cyclic.
   Since $H=RP$ is a $p$-solvable $\mathcal{H}_p$-group with $p$-length 1,
   it follows by Lemma \ref{lem: Hp p-sol gp with p-length 1} that $P\unlhd H$.
   So, $[R,P]=1$,
   and consequently, $[R,G]=1$ by Lemma \ref{lem: G=opG=op'G}.
   So, $R\leq G'\cap \mathbf{Z}(G)$, and hence $R$ is isomorphic to a quotient group of the Schur multiplier $\mathrm{M}(G/R)$.
   Note that $G/R$ is a nonabelian simple $\mathcal{H}_p$-group of order divisible by $p$, 
   and hence Theorem \ref{thm: simple Hp gp} yields that either $G/R\cong \mathrm{PSL}_{2}(p^f)$ with $f\geq 2$ and $p>2$,
   or $G/R\cong \mathrm{PSL}_{3}(4)$ and $p=3$.
   Finally, one checks via \cite{atlas} that either (2) or (3) holds. 
\end{proof}

\begin{lem}\label{lem: under hyp perfect N p}
	Assume Hypothesis \ref{hyp: perfect} and that $p \mid |R|$.
	Set $V=\mathbf{O}_{p}(G)$ and $C=\mathbf{C}_{G}(V)$.
    Then $V\in \mathrm{Syl}_{p}(R)$ is a minimal normal subgroup of $G$, $V\leq C\leq R$, and
	 one of the following holds.  
	 \begin{description}
		\item[(1)] $G$ acts transitively on $V^\sharp$, and we are in one of the three cases below.
		\begin{description}
			\item[(1a)] $p=2$, $(G/V,|V|)=(\mathrm{SL}_2(q),q^2)$ where $q=2^f\geq 4$.
			\item[(1b)] $p>2$, $(G/C,|V|)=(\mathrm{SL}_2(q),q^2)$ where $q=p^f> 4$.
			\item[(1c)] $p=3$, $(G/C,|V|)=(\mathrm{SL}_2(13),3^6)$ and $H^2(G/C,V)=0$.
		\end{description}
	    \item[(2)] $p=3$, $(G/C,|V|)=(\mathrm{SL}_2(5),3^4)$, $H^2(G/C,V)=0$, and the orbit sizes of $G/C$ on $V$ are $1,40,40$.  
	 \end{description}
\end{lem}
\begin{proof}
	By Lemma \ref{lem: under hyp perfect basics},
	we know that $R$ does not have a $G$-chief factor of order $p$.
	As $R$ is the $p$-solvable radical of $G$, $V\leq R$.
	Note that $R$ is also an $\mathcal{H}_p$-group with $p\mid |R|$.
	Thus, an application of Theorem \ref{thm: classification of p-sol Hp} to $R$ yields that $V$ is a nontrivial abelian normal Sylow $p$-subgroup of $R$ with $V\cap \mathbf{Z}(G)=1$.
	As $R$ is the unique maximal normal subgroup of $G$ and $C=\mathbf{C}_{G}(V)\unlhd G$, it follows that $C\leq R$.
	In fact, otherwise $G=C=\mathbf{C}_{G}(V)$, which contradicts $V\cap \mathbf{Z}(G)=1$.
   Let $P\in \mathrm{Syl}_{p}(G)$.
   Since $G=\mathbf{O}^{p'}(G)$ and $V\cap \mathbf{Z}(G)=1$,
   it follows by Lemma \ref{lem: G=opG=op'G} that $[P,V]>1$.
   In particular, $P$ is nonabelian.
   So, by Lemma \ref{lem: basic facts on Hp-group when OpG>1} (3), we conclude that 
   both $V$ and $\mathrm{Irr}(V)$ are irreducible $\mathbb{F}_p[G]$-modules.
   Moreover, $C=\mathbf{C}_{G}(V)=V\times D$ where $D=\mathbf{O}_{p'}(C)$.    
   In particular, $D\unlhd G$.
    Set $U=\mathrm{Irr}(V)$.

   We consider first the case $p=2$.
	As $G/R$ is a nonabelian simple $\mathcal{H}_2$-group,
	it follows by Theorem \ref{thm: simple Hp gp} that $G/R\cong \mathrm{SL}_{2}(q)$ where $q=2^f$ and $f\geq 2$.
     Also, we have $R=V$.
	Indeed, otherwise, 
	 as the $\mathcal{H}_2$-groups $G/V$ and $R/V$ satisfy the hypotheses of Lemma \ref{lem: under hyp perfect N p'},
	 it follows by Lemma \ref{lem: under hyp perfect N p'} that $G/V$ has a cyclic Sylow $2$-subgroup,
	 a contradiction.
    Applying \cite[Lemma 2.5]{qian15}, we deduce that $|V|=q^2$.
	We now claim that $G$ acts transitively on $V^\sharp$.
	To see this, by \cite[Corollary 6.33]{isaacs76}, it suffices to show that $G$ acts transitively on $U^\sharp$.
    Assume that there is some $\lambda\in U^\sharp$ which is fixed by a nontrivial
	subgroup of $G/V$ of odd order $k$.
	  According to \cite[Lemma 2.5]{qian15}, 
	  there exists some $\mu \in U^\sharp$ such that $\mathrm{I}_G(\mu)/V$ contains a Frobenius subgroup of order $2k$. 
	  However, by Lemma \ref{lem: basic facts on Hp-group when OpG>1} (1), $\mathrm{I}_G(\mu)/V$ contains a unique Sylow $2$-subgroup of $G/V$, a contradiction. 
	 Therefore, $\mathrm{I}_{G}(\lambda)\in \mathrm{Syl}_{2}(G)$ for every $\lambda\in U^\sharp$.
	  Equivalently, $G$ acts transitively on $U^\sharp$.
	  Thus, (1a) holds.

	  On the other hand, we assume that $p>2$.
      Recall that $P\in \mathrm{Syl}_{p}(G)$ is nonabelian.
   By Lemma \ref{lem: basic facts on Hp-group when OpG>1}, the action of $G/C$ on $U$ satisfies the hypotheses of Theorem \ref{thm: liebeck 2}.
   If $G$ acts transitively on $U^\sharp$,
   then it also acts transitively on $V^\sharp$ by \cite[Corollary 6.33]{isaacs76}.
   Therefore, (1b) and (1c) follow by Theorem \ref{thm: liebeck 2} and Remark \ref{rmk: smallgroup}.
   Indeed, if $p=3$ and $(G/C,|V|)=(2^{1+4}_-\cdot \mathsf{A}_5,3^4)$,
   as $C=V\times D$, it follows by Remark \ref{rmk: smallgroup} that $G/D$ is not an $\mathcal{H}_3$-group,
   a contradiction.
   Assume now that $G$ does not act transitively on $U^\sharp$.
   We claim that $(G/C,|U|)$ is neither $(M_{11},3^5)$ nor $(\mathrm{PSL}_2(11),3^{5})$.
   In fact, otherwise, $(G/C,|U|) \in  \{ (M_{11},3^5), (\mathrm{PSL}_2(11),3^{5}) \}$;
   in this case, $p=3$ and $R=C=V\times D$ where $D=\mathbf{O}_{3'}(C)$;
   setting $\overline{G}=G/D$, we deduce that $\overline{V}$ is a $5$-dimensional irreducible 
   $\overline{G}/\overline{V}$-module over $\mathbb{F}_3$, where $\overline{G}/\overline{V}\in \{ M_{11}, \mathrm{PSL}_2(11) \} $;
   so $\overline{G}$ is not an $\mathcal{H}_3$-group by Remark \ref{rmk: smallgroup}, a contradiction.
   Applying Theorem \ref{thm: liebeck 2} to $G/C$ and $U$,
   we deduce that $p=3$ and $(G/C,|U|)=(\mathrm{SL}_{2}(5),3^4)$ with orbit sizes $1,40,40$. 
   In this case, $V$ is a $4$-dimensional faithful irreducible $G/C$-module over $\mathbb{F}_3$,
   so (2) holds by Remark \ref{rmk: smallgroup}.
\end{proof}

\begin{lem}\label{lem: under hyp perfect N p C>1}
	Assume Hypothesis \ref{hyp: perfect} and that $p \mid |R|$.
	Set $V=\mathbf{O}_{p}(G)$.
	Then $V=\mathbf{C}_{G}(V)$, and one of the following holds.
	\begin{description}
		\item[(1)] $G/V\cong \mathrm{SL}_2(q)$ where $q=p^f\geq 4$, and $V$ is the natural module for $G/V$.
		\item[(2)] $p=3$, $G=V \rtimes H$ where $H\cong \mathrm{SL}_2(13)$ and $V$ is a $6$-dimensional irreducible $\mathbb{F}_3[H]$-module.
		\item[(3)] $p=3$,
		$G=V \rtimes H$ where $H\cong \mathrm{SL}_2(5)$ and $V$ is a $4$-dimensional irreducible $\mathbb{F}_3[H]$-module.	
	\end{description}
\end{lem}
\begin{proof}
	Assume that $V=\mathbf{C}_{G}(V)$.
	By Lemma \ref{lem: under hyp perfect N p}, Lemma \ref{lem: SL} and Remark \ref{rmk: smallgroup},
	it follows that one of the cases (1), (2) or (3) holds.
	So, it remains to show that $V=\mathbf{C}_{G}(V)$.

	Set $C=\mathbf{C}_{G}(V)$.
	By Lemma \ref{lem: under hyp perfect N p}, we know that $V \in \mathrm{Syl}_{p}(R)$ is minimal normal in $G$ and $V\leq C\leq R$.
	If $V=R$, then we are done.
    So, we may assume that $V<R$.
	Since $C=\mathbf{C}_{G}(V)$, it follows that $C=V\times D$, where $D=\mathbf{O}_{p'}(C)$.
	Given that $G=\mathbf{O}^{p'}(G)$,
	we also have that $D\leq \mathbf{Z}(G)$ by Lemma \ref{lem: basic facts on Hp-group when OpG>1} (2).
	In particular, $[G,C]\leq V$.
	Let $P\in \mathrm{Syl}_{p}(G)$.

	We consider first the case that $PC/C$ is not cyclic.
    By Lemma \ref{lem: under hyp perfect N p'},
	either $G/V\cong \mathrm{SL}_2(p^f)$ with $p>2$ and $f\geq 2$
	or $G/R\cong \mathrm{PSL}_{3}(4)$.
	Applying Lemma \ref{lem: under hyp perfect N p}, 
	we conclude that $G/C\cong \mathrm{SL}_2(p^f)$.
    So, $C=V$.
	
	Next, we assume that $PC/C$ is cyclic.
	By Lemma \ref{lem: under hyp perfect N p}, either $G$ acts transitively on $V^\sharp$ and
	\[
	 (G/C,|V|)\in \{ (\mathrm{SL}_{2}(p),p^2),(\mathrm{SL}_{2}(13),3^{6})\},	
	\]
     or $(G/C,|V|)=(\mathrm{SL}_{2}(5),3^4)$ with orbit sizes $1,40,40$.
     We observe that $|G/C|_p=p>2$ and $G/C\cong \mathrm{SL}_{2}(\ell)$
	 for some odd prime $\ell$ larger than 3.
Therefore, $[G,R]\leq C$.
In particular, $[P,R]\leq C$.
Note that $[P,C]\leq [G,C]\leq V$ and that $P$ acts coprimely on $R/V$, and so $[P,R]\leq V$.
By Lemma \ref{lem: G=opG=op'G}, $R/V\leq \mathbf{Z}(G/V)$.
Since $G/V$ is perfect, $R/V$ is isomorphic to a quotient group of $\mathrm{M}(G/R)$.
Given that $G/R\cong \mathrm{PSL}_{2}(\ell)$ with $\ell$ an odd prime larger than 3,
 we deduce by \cite{atlas} that $|\mathrm{M}(G/R)|=2$.
Therefore,
$C=V$.
\end{proof}

Now, we are ready to classify finite non-$p$-solvable $\mathcal{H}_p$-groups.

\begin{thm}\label{thm: classification of nonsol Hp-gp}
	Let $G$ be a finite non-$p$-solvable
	group.
    Set $N=\mathbf{O}^{p'}(G)$ and $V=\mathbf{O}_{p}(N)$.
	Then $G$ is an $\mathcal{H}_p$-group if and only if one of the following holds.
	\begin{description}
		\item[(1)] $N$ is a nonabelian simple group, and one of the following holds.
		\begin{description}
			\item[(1a)] $p>2$, and $N$ has a cyclic Sylow $p$-subgroup.
			\item[(1b)] $N\cong \mathrm{PSL}_2(q)$ with $q=p^f$ and $f \geq 2$.
			\item[(1c)] $(N,p)\in \{ (\mathrm{PSL}_3(4),3),(M_{11},3),({}^2F_4(2)',5) \}$.
		\end{description}
		\item[(2)] $p>2$, $\mathbf{O}_{p'}(N)>1$, and $N/\mathbf{O}_{p'}(N)$ is a nonabelian simple group, and one of the following holds.
		\begin{description}
			\item[(2a)] $N$ has a cyclic T.I. Sylow $p$-subgroup. 
			\item[(2b)] $N\cong \mathrm{SL}_2(q)$ with $q=p^f$ and $f\geq 2$.
			\item[(2c)] $p=3$, and $N$ is a perfect central extension of $\mathbf{O}_{p'}(N)$ by $N/\mathbf{O}_{p'}(N)\cong \mathrm{PSL}_{3}(4)$.
		\end{description}
		\item[(3)] $V=\mathbf{C}_{N}(V)$, and one of the following holds. 
		\begin{description}
			\item[(3a)] $N/V\cong \mathrm{SL}_2(q)$ where $q=p^f\geq 4$, and $V$ is the natural module for $N/V$.
			\item[(3b)] $p=3$, $N=V \rtimes H$ where $H\cong \mathrm{SL}_2(13)$ and $V$ is a $6$-dimensional irreducible $\mathbb{F}_3[H]$-module.
			\item[(3c)] $p=3$,
			$N=V \rtimes H$ where $H\cong \mathrm{SL}_2(5)$ and $V$ is a $4$-dimensional irreducible $\mathbb{F}_3[H]$-module.	
		\end{description}
	\end{description}
\end{thm}
\begin{proof}
	We first assume that $G$ is a nonsolvable $\mathcal{H}_p$-group.
    Let $M$ be a minimal non-$p$-solvable normal subgroup of $G$.
    We have $M=\mathbf{O}^{p}(M)=\mathbf{O}^{p'}(M)$.
	If $M$ is minimal normal in $G$,
	then case (1) holds by Lemma \ref{lem: nonsol Hp gp with nonabel minimal normal} and Theorem \ref{thm: simple Hp gp}.
	Assume now that $M$ is not minimal normal in $G$.
    Let $R$ be the $p$-solvable radical of $M$.
	Then $G/R$ satisfies the hypotheses of Lemma \ref{lem: nonsol Hp gp with nonabel minimal normal},
	and so $M/R=\mathbf{O}^{p'}(G/R)$ is a nonabelian simple $\mathcal{H}_p$-group of order divisible by $p$.
	Hence, $N=\mathbf{O}^{p'}(G)=M$ satisfies Hypothesis \ref{hyp: perfect}.
	Applying Lemmas \ref{lem: under hyp perfect N p'} and \ref{lem: under hyp perfect N p C>1} to $N$,
	we conclude that either (2) or (3) holds.

	Conversely, suppose that one of the cases (1), (2) or (3) holds.
    Let $\chi \in \mathrm{Irr}(G)$ be of degree divisible by $p$.
	To see that $G$ is an $\mathcal{H}_p$-group, it is enough to show that $\mathrm{cod}(\chi)_p=1$.
	Let $\theta$ be an irreducible constituent of $\chi_N$.
	As $N=\mathbf{O}^{p'}(G)$, it follows that  $\chi(1)_p=\theta(1)_p$.

    Assume that  either (1) or (2) holds.
    Then $N/\mathbf{O}_{p'}(N)$ is a nonabelian simple $\mathcal{H}_p$-group by Theorem \ref{thm: simple Hp gp}.
    Therefore, $N/\mathbf{O}_{p'}(N)$ is an $\mathcal{H}^*_p$-group.
    Now, we claim that $N$ is also an $\mathcal{H}^*_p$-group.
	Indeed, if (1) holds, then we are done; 
	if (2a) holds, then $N$ is an $\mathcal{H}^*_p$-group by Proposition \ref{prop: tisylow};
	if either (2b) or (2c) holds,
	as $\mathbf{O}_{p'}(N)\leq \mathbf{Z}(N)$, then $N$ is an $\mathcal{H}^*_p$-group by Corollary \ref{cor: subgroup is Hp*} (2).
	Noting that $N=\mathbf{O}^{p'}(G)$ is an $\mathcal{H}^*_p$-group,
	we conclude by Corollary \ref{cor: subgroup is Hp*} (3) that 
	$G$ is an $\mathcal{H}^*_p$-group.

    Assume that (3) holds. 
    We assert that $p\nmid \varphi(1)$ for each $\varphi \in \mathrm{Irr}(N|V)$.
    In fact, this assertion is verified by \cite[Proposition 2.3]{gold2014} when (3a) holds,
	and by $\mathsf{GAP}$ \cite{gap} when either (3b) or (3c) holds.
	Consequently, we have $V\leq \ker(\theta)$.
    Since $\chi\in \mathrm{Irr}(G|\theta)$, it follows that $V\leq \ker(\chi)$, i.e. $\chi \in \mathrm{Irr}(G/V)$.
    	Noting that either $N/V\cong \mathrm{SL}_2(q)$ where $q=p^f\geq 4$ or $|N/V|_p=p$,
		we deduce that $N/V$ has an abelian T.I. Sylow $p$-subgroup, and therefore
		$N/V$ is an $\mathcal{H}^*_p$-group by Proposition \ref{prop: tisylow}.
	   Observing that $N/V=\mathbf{O}^{p'}(G/V)$, we conclude that $G/V$ is also an $\mathcal{H}^*_p$-group
		by Corollary \ref{cor: subgroup is Hp*} (3). 
   Therefore, $\chi(1)_p=|G/V|_p$.
	In particular, $\mathrm{cod}(\chi)_p=1$.
	\end{proof}

Finally, we give a rough description of the groups arising in the subcase (2a) of Theorem \ref{thm: classification of nonsol Hp-gp}.

\begin{thm}\label{thm: (6a) of thmA}
	Let $G=\mathbf{O}^{p'}(G)$ be a finite group where $p$ is an odd prime and let $P\in \mathrm{Syl}_{p}(G)$.
	Assume that $G/\mathbf{O}_{p'}(G)$ is a nonabelian simple group.
	Then $G$ has a cyclic T.I. Sylow $p$-subgroup if and only if one of the following holds.
	\begin{description}
		 \item[(1)] The $p$-solvable group $H:=\mathbf{O}^{p'}(\mathbf{O}_{p'}(G) P)=H' \rtimes P$
		 where $P$ is a cyclic T.I. subgroup of $H$. 
		 \item[(2)] $G$ is a quasisimple group with a cyclic Sylow $p$-subgroup.
	\end{description}
\end{thm} 
\begin{proof}
    We first assume that $G$ has a cyclic T.I. Sylow $p$-subgroup.
	Then $G$ is an $\mathcal{H}^*_p$-group by Proposition \ref{prop: tisylow}.
	If $[P,\mathbf{O}_{p'}(G)]=1$, as $G=\mathbf{O}^{p'}(G)$ and $G/\mathbf{O}_{p'}(G)$ is nonabelian simple,
	then $\mathbf{O}_{p'}(G)\leq \mathbf{Z}(G)\cap G'$.
	Therefore, case (2) holds.
	Assume now that $[P,\mathbf{O}_{p'}(G)]>1$.
    Since $\mathbf{O}_{p'}(G) P$ is an $\mathcal{H}^*_p$-group with $p$-length 1 by Corollary \ref{cor: subgroup is Hp*} (1),
	it follows that case (1) holds by Lemma \ref{lem: Hp p-sol gp with p-length 1}.

	We assume next that either (1) or (2) holds.
	Set $\overline{G}=G/\mathbf{O}_{p'}(G)$.
    Then the nonabelian simple group $\overline{G}$ has a cyclic Sylow $p$-subgroup $\overline{P}$.
	By \cite[Theorem 1]{blau85},
	we know that $\overline{P}$ is a T.I. subgroup of $\overline{G}$.
    Let $x\in G$.
	Then either $\overline{P^x}\cap \overline{P}=1$ or $\overline{P^x}=\overline{P}$.
	If the former holds, then $P^x\cap P\leq \mathbf{O}_{p'}(G)\cap P=1$.
		Assume now that the latter holds.
	Then $\mathbf{O}_{p'}(G) P^{x} =\mathbf{O}_{p'}(G)P$.
	If case (1) holds, as $P$ and $P^x$ are T.I. Sylow $p$-subgroups of $H=\mathbf{O}^{p'}(\mathbf{O}_{p'}(G) P)$, 
	then either $P^x\cap P=1$ or $P^x=P$.
	If case (2) holds, as $ \mathbf{O}_{p'}(G) P= \mathbf{O}_{p'}(G) \times P$,
	then $P^x=P$.
	Consequently, $P$ is a cyclic T.I. Sylow $p$-subgroup of $G$.
\end{proof}

\begin{acknowledgement}
	The authors would like to thank Professor Silvio Dolfi for some 
	suggestions on English writing.
	The authors are grateful to  the referee for her/his
 valuable comments.
\end{acknowledgement}



\begin{thebibliography}{ABCDEF}\setlength{\itemsep}{-2mm} 
\small

	

	\bibitem[APS24]{akhlaghi24}
	Z. Akhlaghi, E. Pacifici and L. Sanus,
	\newblock Element orders and codegrees of characters in non-solvable groups,
	\newblock \emph{J. Algebra} {\bf 644} (2024), 428--441.

	\bibitem[Bla85]{blau85}
	H.I. Blau,
	\newblock On trivial intersection of cyclic Sylow subgroups,
	\newblock \emph{Proc. Amer. Math. Soc.} {\bf 94(4)} (1985), 572--576.
  
  
	\bibitem[BM90]{blau90}
	H.I. Blau and  G.O. Michler,
	\newblock Modular representation theory of finite groups with T.I. Sylow $p$-subgroups,
	\newblock \emph{Trans. Amer. Math. Soc.} {\bf 319(2)} (1990), 417--468.




	\bibitem[Car72]{carter72}
	R.W. Carter,
	\newblock \emph{Simple Groups of Lie type},
	\newblock John Wiley \& Sons, London, 1972.

	
	\bibitem[CN22]{chen22}
	X. Chen and G. Navarro,
	\newblock Brauer characters, degrees and subgroups,
	\newblock \emph{Bull. London Math. Soc.} {\bf 54(3)} (2022), 891--893.


	\bibitem[CH89]{chillag89}
	D. Chillag and M. Herzog,
	\newblock On character degrees quotients,
	\newblock \emph{Arch. Math.} {\bf 55(1)} (1989), 25--29.


	\bibitem[CMM91]{chillag91}
	D. Chillag, A. Mann, and O. Manz,
	\newblock The co-degrees of irreducible characters,
	\newblock \emph{Israel J. Math.} {\bf 73(2)} (1991), 207--223.


	\bibitem[Atl1]{atlas}
	J.H. Conway, R.T. Curtis, S.P. Norton, R.A. Parker and R.A. Wilson,
	\newblock {\em $\mathbb{ATLAS}$ of finite groups: maximal subgroups and ordinary characters for simple groups},
	\newblock Clarendon Press, Oxford, 1985.






	








	\bibitem[GL99]{gagola99}
	S.M. Gagola and  M.L. Lewis,
	\newblock A character theoretic
	condition characterizing nilpotent groups,
	\newblock \emph{Comm. Algebra} {\bf 27(3)} (1999), 1053--1056.


	
	\bibitem[GAP]{gap} The GAP Group, GAP - Groups, Algorithms, and Programming, Version 4.13.1, 2024, \href{http://www.gap-system.org}{http://www.gap-system.org}.

	\bibitem[Gia24]{giannelli24}
	E. Giannelli,
	\newblock Character codegrees and element orders in symmetric and alternating groups,
	\newblock \emph{J. Algebra Appl.}, {\bf 23(09)} (2024), 2450144.
	
	



	\bibitem[GLP+16]{giudici16}
	M. Giudici, M.W. Liebeck, C.E. Praeger, J. Saxl and P.H. Tiep,
	\newblock Arithmetic results on orbits of linear groups,
	\newblock \emph{Trans. Amer. Math. Soc.} {\bf 368(4)} (2016), 2415--2467.


	\bibitem[GGL+14]{gold2014}
	D. Goldstein, R. Guralnick, M. Lewis, A. Moret\'{o}, G. Navarro and P.H. Tiep,
	\newblock Groups with exactly one irreducible character of degree divisible by $p$,
	\newblock \emph{Algebra Number Theory} {\bf 8(2)} (2014), 397--428.


	\bibitem[GLS94]{gorenstein94}
	D. Gorenstein, R. Lyons and R. Solomon,
	\newblock {\em The classification of the finite simple groups, number 3},
	\newblock American Mathematical Society, 1994.
	



	
	\bibitem[HM25]{hung25}
	N.N. Hung and A. Moret\'o,
	\newblock The codegree isomorphism problem for finite simple groups II,
	\newblock \emph{Quart. J. Math.} {\bf 76(1)} (2025), 237--250.
	
	
	

	\bibitem[Hup67]{huppert67}
	B. Huppert,
	\emph{Endliche gruppen I}, Springer Verlag, Berlin, 1967.




	\bibitem[Isa76]{isaacs76}
	I.M. Isaacs,
	\emph{Character Theory of Finite Groups}, Academic Press, New York, 1976.


	\bibitem[Isa11]{isaacs11}
	I.M. Isaacs,
	\newblock Element orders and character codegrees,
	\newblock \emph{Arch. Math.} {\bf 97(6)} (2011), 499--501.




	\bibitem[KM13]{kessar13}
	R. Kessar and G. Malle,
 \newblock Quasi-isolated blocks and Brauer's height zero conjecture, 
 \newblock \emph{Ann. Math.} {\bf 178(1)} (2013), 321--384. 



 \bibitem[KM25]{kourovka notebook}
 E.I. Khukhro and V.D. Mazurov,
 \newblock {\em The Kourovka notebook: Unsolved problems in group theory}, 20th edition,
 \newblock \href{https://kourovka-notebook.org}{https://kourovka-notebook.org}.



	
	
	 \bibitem[LQ16]{liang16}
	 D. Liang and G. Qian,
	 \newblock Finite groups with coprime character degrees and codegrees,
	 \newblock \emph{J. Group Theory} {\bf 19(5)} (2016), 763--776.
 

	 \bibitem[Lie87]{liebeck87}
	 M.W. Liebeck,
	 \newblock The affine permutation groups of rank three,
	 \newblock \emph{Proc. London Math. Soc.} {\bf 54(3)} (1987), 477--516.


	\bibitem[Liu22]{liu22}
	Y. Liu,
	\newblock Nonsolvable groups whose irreducible character degrees have
	special $2$-parts,
	\newblock \emph{Front. Math. China} {\bf 17(6)} (2022), 1083--1088.



	\bibitem[Mad23]{madanha23}
	S.Y. Madanha,
	\newblock Codegrees and element orders of almost simple groups,
	\newblock \emph{Comm. Algebra} {\bf 51(7)} (2023), 3143--3151.





	\bibitem[MN21]{malle21}
	G. Malle and G. Navarro,
	\newblock Brauer's height zero conjecture for principal blocks,
	\newblock \emph{J. reine angew. Math.} {\bf 778} (2021), 119--125.


	



	\bibitem[MW92]{manzwolfbook}
	O. Manz and T.R. Wolf,
	\newblock {\em Representations of Solvable Groups},
	\newblock Cambridge University Press, Cambridge, 1992.

	
	

	
	
	\bibitem[Mic86]{michler86}
	G. Michler,
	\newblock A finite simple Lie group has $p$-blocks with  different defects,
	$p\neq 2$,
	\newblock \emph{J. Algebra} {\bf 104(2)} (1986), 220--230.

	\bibitem[Mor05]{moreto05}
	A. Moret\'{o},
	\newblock Complex group algebras of finite groups: Brauer's problem 1,
	\newblock \emph{Electron. Res. Announc. Amer. Math. Soc.} {\bf 11} (2005), 34--39.




  \bibitem[MH24]{moreto24}
	A. Moret\'{o} and N.N. Hung,
	\newblock The codegree isomorphism problem for finite simple groups,
	\newblock \emph{Quart. J. Math.} {\bf 75(3)} (2024), 1157--1179.
	
	
	
	

 
 

	\bibitem[Nav98]{navarrobook}
	G. Navarro,
	\emph{Characters and blocks of finite groups}, Cambridge University Press, Cambridge, 1998.




	\bibitem[P\'al01]{palfy01}
	P. P\'alfy,
	\newblock On the character degree of solvable groups, II: disconnected graphs,
	\newblock \emph{Studia Sci. Math. Hungar.} {\bf 38} (2001), 339--355.


	\bibitem[PR02]{parkerbook}
	C. Parker and P. Rowley,
	\emph{Symplectic Amalgams},  Springer-Verlag London, Ltd., London, 2002.




	\bibitem[PT91]{pazderski91}
	G. Pazderski and G. Tiedt,
	\newblock On groups with extremal $p$-blocks,
	\newblock \emph{Arch. Math.} {\bf 57(3)} (1991), 216--220.




	\bibitem[Qia02]{qian02}
	G. Qian,
	\newblock Notes on character degrees quotients of finite groups,
	\newblock \emph{J. Math.(in Chinese)} {\bf 22(2)} (2002), 217--220.

	\bibitem[Qia11]{qian11}
	G. Qian,
	\newblock A note on element orders and character codegrees,
	\newblock \emph{Arch. Math.} {\bf 97(2)} (2011), 93--103.

	

	

	\bibitem[Qia12]{qian12}
	G. Qian,
	\newblock A character theoretic criterion for $p$-closed group,
	\newblock \emph{Israel J. Math.} {\bf 190(1)} (2012), 401--412.


	\bibitem[Qia21]{qian21}
	G. Qian,
	\newblock Element orders and character codegrees,
	\newblock \emph{Bull. London Math. Soc.} {\bf 53(3)} (2021), 820--824.


	\bibitem[QS04]{qian04}
	G. Qian and W. Shi,
	\newblock The largest character degree and the Sylow subgroups of finite groups, 
	\newblock \emph{J. Algebra} {\bf 277(1)} (2004), 165--171.   

	\bibitem[QWW07]{qian07}
	G. Qian, Y. Wang and H. Wei,
	\newblock Co-degrees of irreducible characters in finite groups,
	\newblock \emph{J. Algebra} {\bf 312(2)} (2007), 946--955.

	


	\bibitem[QY15]{qian15}
	G. Qian and Y. Yang,
	\newblock Nonsolvable groups with no prime dividing three character degrees,
	\newblock \emph{J. Algebra} {\bf 436} (2015), 145--160.




	\bibitem[Sch85]{schmid85}
	P. Schmid,
	\newblock Rational matrix groups of a special type,
	\newblock \emph{Linear Algebra Appl.} {\bf 71} (1985), 289--293.

	\bibitem[Ton25]{tongviet25}
	H.P. Tong-Viet,
	\newblock A characterization of some finite simple groups by their character codegrees,
	\newblock \emph{Math. Nachr.} {\bf 298(4)} (2025), 1356--1369.

	


	\bibitem[Zha00]{zhang00}
	J. Zhang,
	\newblock A note on character degrees of finite solvable groups,
	\newblock \emph{Comm. Algebra} {\bf 28(9)} (2000), 4249--4258.




\end{thebibliography}
\end{document}